\documentclass[12pt]{article}
\usepackage{amsmath, amssymb,epsfig,bm}
\usepackage{latexsym}
\usepackage{algorithm,algpseudocode}
\usepackage{graphicx}
\usepackage{subfigure}
\usepackage{color}
\usepackage{hyperref}
\numberwithin{equation}{section}
\def\debproof{ {\bf Proof.} }
\def\finproof{\hfill {\small $\Box$} \\}

\newcommand{\bE}{{\bf E}}

\newcommand{\bJ}{{\bf J}}
\newcommand{\bF}{{\bf F}}

\newcommand{\bn}{{\bf n}}

\newcommand{\bp}{{\bf \Phi_\alpha}}

\newtheorem{lem}{Lemma}[section]

\newtheorem{thm}{Theorem}[section]

\begin{document}	

\title{An iterative Method for the Inverse Eddy Current Problem with Total Variation Regularization}

\author{Junqing Chen\thanks{\footnotesize
Department of Mathematical Sciences, Tsinghua University, Beijing
100084, China. The work of this author was partially supported by NSFC under the grant 11871300 and by National Key  R\&D Program of China 2019YFA0709600, 2019YFA0709602. (jqchen@tsinghua.edu.cn).}
\and Zehao Long\thanks{\footnotesize Department of Mathematical Sciences, Tsinghua University, Beijing 100084, China. 
(longzh18@mails.tsinghua.edu.cn).}
}
\date{}
\maketitle

\begin{abstract}
Conductivity reconstruction in an inverse eddy current problem is considered in the present paper. With the electric field measurement on part of domain boundary, we formulate the reconstruction problem to a constrained optimization problem with total variation regularization. Existence and stability are proved for the solution to the optimization problem. The finite element method is employed to discretize the optimization problem. The gradient Lipschitz properties of the objective functional are established for the  the discrete optimization problems. We propose a modified alternating direction method of multipliers  to solve the discrete problem and prove the convergence of the algorithm. Finally, we show some numerical experiments to illustrate the efficiency of the proposed method.
\end{abstract}
			


\section{Introduction}\label{sect1}
Eddy current inversion is an important non-destructive testing modality and  has been used in a wide range of applications such as geophysical prospecting, flaw detection and safety inspection \cite{ABF,ACCVW,HHT}. This inversion method uses induced electromagnetic data to detect the conducting anomalies or flaw in conductive objects. Comparing with the inversions using acoustic wave or elastic wave,  the eddy current method is more sensitive to conductive medium.  The induced electromagnetic field is modeled by Maxwell's equations in low frequency. While in low frequency case, the eddy current model is a good approximation \cite{BAN} and can be solved with high efficient algorithm such as \cite{CCCZ,HX}.
	
There are tremendous efforts devoted to solve the inverse problems related to eddy current.  Here we list some literature for the relation work. An inverse source problem for eddy current equation has been discussed in \cite{Rodr2012Inverse}. To detect and recognize small conducting anomalies, the conductive Generalized Polarization Tensors (GPTs) for small inclusion are studied and a MUSIC-like algorithm is proposed  in \cite{ACCGV,ACCVW}.  In \cite{CLZ},  mathematical and numerical theory of an inverse eddy current problem has been studied and an NLCG algorithm has been proposed to reconstruct the conductivity. A monotone method is introduced in \cite{TG} to recover the conductivity.  
Meanwhile, there are many works about the inverse eddy problem in industry, too. Specially in geophysical prospecting field, we refer to the monographs \cite{Haber2014Computational} and \cite{Zhdanov}, to name a few. 
 
\begin{figure}
	\centering{
		\includegraphics[width=0.4\textwidth]{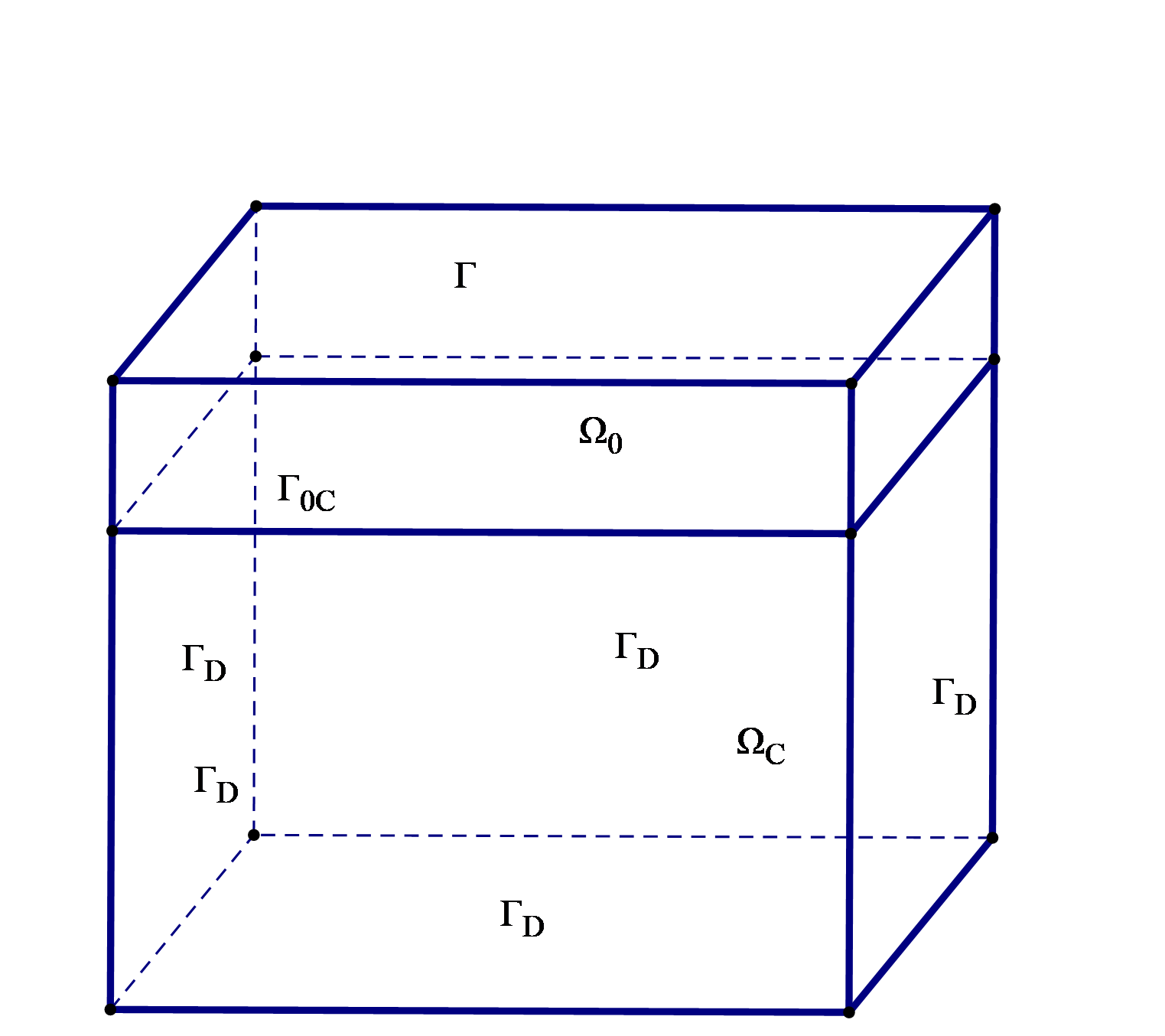}
	}
	\caption{geometric setting}\label{fig:1}
\end{figure}

 To specify the problem, we  depict in Figure \ref{fig:1} the typical domain $\Omega$ for the eddy current problem.  The domain $\Omega$ contains three parts:  $\Omega_0$, $\Omega_1$ and $\Omega_2$, where $\Omega_0$ is non-conducting part, $\Omega_1,\Omega_2$ are conducting parts.  
The forward model is governed by the following eddy current equations
 \begin{equation}\label{eddy}
	\left\{
	\begin{aligned}
		\nabla\times(\mu^{-1}\nabla\times \bE(\sigma))-i\omega(\sigma+\sigma_0) \bE(\sigma)=&i\omega \bJ_s &\mbox{ in }\Omega,\\
		\nabla\cdot \epsilon \bE(\sigma)=&0 &\mbox{ in }{\Omega_0},
	\end{aligned}
	\right.
\end{equation}    
with the mixed boundary conditions 
\begin{equation}\label{bdry}
	\bn\times\nabla\times \bE(\sigma)=0, \quad \bn\cdot \bE(\sigma)=0\mbox{ on }\Gamma; \quad \bn\times \bE(\sigma)=0\mbox{ on }\Gamma_D. 
\end{equation}
Where $\bJ_s$ is the source satisfying $\nabla\cdot (\bJ_s(x))=0\quad x\in \Omega$ and $\bn$ is the unit outer normal vector of $\Gamma$.  $\Gamma$ is the upper boundary of $\Omega$ and $\Gamma_D$ is the rest part of the boundary of $\Omega$. 
 We assume that magnetic permeability $\mu$ is a constant function in $\Omega$, $\sigma_0$ is the background conductivity which vanishes in $\Omega_0$ and is a positive constant in $\Omega_C=\Omega_1\cup\overline{\Omega}_2$ and $\sigma$ is the unknown conductivity which is supported in $\Omega_2$. This means that $\sigma$ is compactly supported in $\Omega_C$. The real conductivity $\sigma_0+\sigma$ has positive lower bound in $\Omega_C$.  $\varepsilon$ is the electric permittivity in $\Omega_0$.

As for the inverse problem, with the data $\bn\times\bE$ on $\Gamma$, we want to reconstruct the unknown conductivity $\sigma$. In \cite{CLZ},  this problem is modeled as a constrained optimization problem with eddy current equations as constraint. With the assumption that $\sigma\in H^1_0(\Omega_C)$, they study the well-posedness of the inverse problem. To improve the stability of the inverse problem, the $H^1$-regularization is introduced.  Usually, the unknown conductivity $\sigma$ is not smooth and does not lie in $H^1_0(\Omega_C)$. In this situation, the $H^1$-regularization is not suitable.  For reconstructing non-smooth parameter, it is known that total variation regularization is a good choice. Comparing with $H^1$ regularization, total regularization can deal with non-smooth parameters  theoretically and , as a consequence, can bring clearer boundary of conductors  numerically. Here we give some reference about the total-variation regularization, such as\cite{ chenzou,GJL,ROF2,ROF1}. 
 The total variation of  a function is defined as 
 $$|Du|(\Omega) = \sup\{\int_{\Omega} u \nabla \cdot g dx : g\in (C_{0}^1(\Omega))^3\quad\text{and}\quad |g(x)|\leq 1 \quad\text{in}\quad \Omega \},$$
which is not differentiable as $L^2$ and $H^1$  regularization terms. This brings us difficulty in solving the optimization problem by iterative method based on gradient.  One way to deal with this difficulty is replacing $TV$ term $|Du|$ by its smoothness $\sqrt{|\nabla u|^2+\varepsilon^2}$ for small $\varepsilon > 0$\cite{chenzou, ROF1}. Another way is employing the splitting method \cite{splitting2,GO,ROF2}. 

To deal with the discontinuity of the parameter $\sigma$, we introduce the total variation regularization to the inverse eddy current problem. In the present work, we treat the $TV$ regularization with splitting method. Then we use the alternating direction method of multiplies (ADMM) to solve the regularized problem. It is known that the ADMM algorithm has been successfully applied into a range of optimization problems such as signal processing, networking,  machine learning problems and so on \cite{ADMM1,ADMM3,ADMM2}. The convergence of ADMM algorithm can be proved under very mild conditions\cite{ADMM1}. Recent works on its rate of convergence can be found in \cite{HY}. Specially, for the nonconvex optimization problem, the key ingredient of convergence analysis is the gradient Lipschitz continuous condition \cite{HLR}. As for the optimization problems concerned inversion,  to prove this condition is a difficult task since the objective function depends on the parameters through partial differential equation. A Lipschitz-like property of gradient can be proved in the inverse discrete eddy current problem , which guarantees the convergence for our algorithm. 

The main contributions are as follows. Firstly, we introduce the TV regularization to the inverse eddy current problem and show the existence and stability of solution to the optimization problem and the corresponding finite element discretization problem. Secondly, we prove that the discrete objective functional has the gradient Lipschitz-like property. Then we introduce an algorithm of ADMM form, which is different from the traditional ADMM method, and prove the convergence for the algorithm. Finally,  we give some numerical experiments to verify the feasibility and convergence of the algorithms.

The outline of this paper is as follows. In Section \ref{sect2}, some theoretical analysis of the inverse problem is presented. We consider the existence and stability of the optimization problem.  In  Section \ref{sect3}, the discrete inverse problem is considered,  the eddy current model, the conductivity, and total variation semi-norm regularization is discretized with finite element method. Then the discrete optimization problem is introduced and the existence and stability of solution is analyzed. Also we prove some inequalities about the objective functional. In Section \ref{sect4}, we propose the modified ADMM algorithm for the discrete inverse problem and prove the convergence. In Section \ref{sect5}, we show some numerical examples to illustrate the effectiveness of proposed algorithm. Finally, in Section \ref{sect6}, we draw some conclusions on this work.

\section{Preliminary results for the inverse eddy current problem}\label{sect2}

We first introduce some Sobolev spaces used in formulating the inverse eddy current problem.
$$
\begin{aligned}
	H_{\Gamma}(\operatorname{curl} ; \Omega) &=\left\{\mathbf{u} \in L^{2}(\Omega)^{3} \mid \nabla \times \mathbf{u} \in L^{2}(\Omega)^{3}, {n} \times \mathbf{u}=0 \text { on } \Gamma_{D}\right\}, \\
	H_{\Gamma}^{1}\left({\Omega_0}\right) &=\left\{v \in L^{2}\left({\Omega_0}\right)\left|\nabla v \in L^{2}\left(\Omega_{0}\right)^{3}, v\right|_{\partial {\Omega_0} \backslash \Gamma}=0\right\}, \\
	{Y} &=\left\{\mathbf{u} \in H_{\Gamma}(\operatorname{curl} ; \Omega) \mid(\varepsilon \mathbf{u}, \nabla \phi)=0 \quad \forall \phi \in H_{\Gamma}^{1}\left(\Omega_0\right)\right\},
\end{aligned}
$$
We define the sesquilinear form $a: H_{\Gamma}(\operatorname{curl}; \Omega)\times H_{\Gamma}(\operatorname{curl}; \Omega) \to \mathbb{C}$ as 
\begin{equation}\label{se-form}
a(\bE, \bF)=\int_{\Omega} \mu^{-1} \nabla \times \bE \cdot \nabla \times \overline{\bF}-i \omega(\sigma+\sigma_0) \bE \cdot \overline{\bF} d x \quad \forall \bE, \bF \in H_{\Gamma}(\operatorname{curl} ; \Omega).
\end{equation}
Then the weak form of equations \eqref{eddy}-\eqref{bdry} is : Find $\bE \in {Y}$, such that
\begin{equation}\label{seq}
	a(\bE, \bF)=i \omega \int_{\Omega} \bJ_s \cdot \overline{\bF} d x \quad \forall \bF \in {Y}.
\end{equation}
The equivalent saddle-point problem to \eqref{seq} is: Find $\bE \in  H_\Gamma(\operatorname{curl};\Omega), \phi \in H^1_\Gamma(\Omega_0)$, such that
\begin{equation}\label{weakeddy}
\left\{\begin{aligned}
		a(\bE, \bF) + b(\nabla \phi, \bF) dx&=i \omega \int_{\Omega} \bJ_s \cdot \overline{\bF} d x \quad &&\forall \bF \in H_\Gamma(\operatorname{curl};\Omega),\\
		b(\bE ,\nabla \psi)&=0 \quad  &&\forall \psi\in H^1_\Gamma({\Omega_0}),
\end{aligned}\right.
\end{equation}
where $b(\bE,\bF) = \int_{\Omega_0} \varepsilon \bE\cdot\overline{\bF} dx,\,\,\forall\,\, \bE,\bF\in H_{\Gamma}(\operatorname{curl},\Omega)$.

We collect some results in \cite{CLZ} to the following theorem. It gives the theoretical supports of the forward problem and also shows the regularity of the solution to the inverse eddy current problem. With this theorem, the data has at least $L^2$ regularity on the boundary, which makes the objective functional in the following parts reasonable. 
\begin{thm} \label{base}
If $(\sigma+\sigma_0)|_{\Omega_C}$ has positive lower bound, then we have the following results.
\begin{itemize}
    \item \textbf{well-posedness of problem (\ref{seq}):} There is a unique $\bE \in {Y}$, such that 
	\begin{equation*}
		a(\bE, \bF)=i \omega \int_{\Omega} \bJ_s \cdot \overline{\bF} d x \quad \forall \bF \in {Y},
	\end{equation*}
	and
	$$\|\bE\|_{H(\operatorname{curl},\Omega)}\leq C\|\bJ_s\|_{L^2(\Omega)}.$$
	 for some const $C>0$, which only depend on $\Omega$ and the lower bound of $(\sigma+\sigma_0)|_{\Omega_C}$.
    \item \textbf{regularity of data:} Assuming that $\Omega_0$ and $\Omega_c$ are polyhedral domains and $\Omega$ is convex, $\sigma_0$ is a constant in $\Omega_C$ and $\mathbf{E}^{o b s}$ is the solution to with the exact conductivity $\sigma_0+\sigma_e$, then for any $\sigma$ we have
   $$
    \left.\left(\mathbf{E}(\sigma)-\mathbf{E}^{o b s}\right)\right|_{\Omega_0} \in {H}^{1 / 2}\left(\Omega_0\right) .
   $$
\end{itemize}
\end{thm}

With the measured electric field on $\Gamma$, we want to  reconstruct $\sigma$ in $\Omega$. This reconstruction, or in other word inversion, can be achieved by minimizing the functional
\begin{equation}
	G(\sigma)=\frac{1}{2}\|\bE(\sigma)\times n-\bE^{obs}\times n\|_{L^2(\Gamma)}^2, \label{obj_G}
\end{equation}
where $\bE(\sigma)$ is the solution of (2.1) and (2.2) for given $\sigma$ and $\bE^{obs}$ denotes measurement on $\Gamma$. 

It is known that the inversion using measurements on $\Gamma$ is an ill-posed problem. 
Then $\|\nabla\sigma\|_{L^2(\Omega_C)}^2$ regularization is introduced to the objective function \eqref{obj_G} and the following constrained optimization problem is proposed to solve the inversion problem \cite{CLZ}.
\begin{eqnarray}\label{H1opt}
&& ~~	 \min_{\sigma\in H^1_0(\Omega_C)}G(\sigma) + \alpha \|\nabla\sigma\|_{L^2(\Omega_C)}^2,\label{pre-problem}\\
&&  \mbox{ s.t. } \bE(\sigma) \mbox{ solves \eqref{eddy}-\eqref{bdry} or \eqref{seq} in weak sense}.\nonumber
\end{eqnarray}

This methods has abundant theorical analysis and successful performance in recovering smooth parameters. But as we mentioned in the introduction, the parameter $\sigma$ is not smooth in many cases,  even when $\sigma$ is piecewise constant. Then in order to  deal with such situations and achieve sharp inversion result with discontinuity,  instead of \eqref{pre-problem},  we employ the BV semi-norm regularization (namely, TV-regularization) and  introduce the following optimization problem to reconstruct the parameter $\sigma$. 
\begin{equation}\label{myproblem0}
	\mathop{\min}_{\sigma \in BV(\Omega_C)}
	\bp(\sigma)= G(\sigma)+\alpha |D\sigma|(\Omega_C), 
\end{equation}
where $\alpha > 0$ is the regularization parameter and $BV(\Omega_C)$ is the space of function with bounded variation, i.e.
$$BV(\Omega_C)=\{u\mid \|u\|_{BV}:=\|u\|_{L^1(\Omega_C)}+|Du|(\Omega_C)<+\infty\}.$$
Let $M>m>-\sigma_0$ and 
$$
K=\{\sigma| \sigma \mbox{ is compactly supported in }\Omega_C\mbox{ and }  m\leq \sigma \leq M \}\cap BV(\Omega_C).
$$
It is clear that with $\sigma\in K$ and the same assumptions on domain $\Omega,\Omega_C$, the results in Theorem \ref{base} are still right. 

To solve problem \eqref{pre-problem} with iterative method,  
the gradient of $G(\sigma)$ is very important.
For any $\sigma\in K$, let $\tilde{\sigma}\in BV(\Omega_C)\cap L^\infty(\Omega_C)$ such that  $\sigma+t\tilde{\sigma}\in K$ for some $t>0$, we can define the G\^ateaux derivative of $G$ at $\sigma$ alone direction $\tilde{\sigma}$ as

\begin{equation}
    G'(\sigma;\tilde{\sigma})=\lim_{t\downarrow 0}(G(\sigma+t\tilde{\sigma})-G(\sigma))/t.\label{gateaux}
\end{equation}
If $\sigma$ is an interior point of $K$,  $\tilde{\sigma}$ can be any element in $BV(\Omega_C)\cap L^\infty(\Omega_C)$ since $K$ is convex.

With the definition of $G(\sigma)$,  the adjoint-state equation to \eqref{weakeddy} is
\begin{equation}\label{adjoint}
    \left\{
    \begin{aligned}       a(\bF,\tilde{\bE})+b(\nabla\psi,\tilde{\bE}) =& \int_{\Gamma} \bn\times \overline{(\bE^{obs} - \bE)}\times \bn \cdot \overline{\tilde{\bE}}ds \qquad &&\forall \tilde{E} \in H_{\Gamma}(\operatorname{curl},\Omega),\\
        b(\bF,\nabla\tilde{\phi})=&0 &&\forall \tilde{\psi} \in H_{\Gamma}^1(\Omega_0),
    \end{aligned}
    \right.
\end{equation}
where $\bE^{obs}$ is the data on $\Gamma$, $\bE(\sigma)$ is the solution to \eqref{weakeddy}.

 Authors in \cite{CLZ} give a formulation of the G\^ateaux derivative \eqref{gateaux} when $\sigma$ is a solution to \eqref{H1opt} through Lagrangian. In fact we can prove the formulation is vaild for all $\sigma \in K$ by definition.
With the help of adjoint state equation, the following lemma gives an explicit representation of the G\^ateaux derivative. Actually, it gives the gradient of $G(\sigma)$ in the weak sense when $\sigma$ is an interior point of $K$.
\begin{lem}\label{derivative}
  Assume $\sigma \in K$, $\tilde{\sigma}\in BV(\Omega_C)\cap L^\infty(\Omega_C)$ such that $\sigma+t\tilde{\sigma}\in K$ for some $t>0$, 
  $\bE(\sigma)$ is the solution to (\ref{weakeddy}) and  $(\bF(\sigma),\psi(\sigma))$ is the  solution to (\ref{adjoint}). Then 
  $$G'(\sigma;\tilde{\sigma}) = \omega\int_{\Omega_C}\mathbf{Im}(\bE\cdot \bF)\tilde{\sigma}dx.$$
\end{lem}
\debproof
For the sake of simplicity, we let $\hat\sigma =\sigma+t\tilde{\sigma}$. Let
$(\bE(\sigma),\phi(\sigma))$ and  $(\bE(\hat\sigma),\phi(\hat\sigma))$ be the solutions to (\ref{weakeddy}) with respect to $\sigma$ and $\hat\sigma$. 
Let $(\bF(\sigma),\psi(\sigma))$ and $(\bF(\hat\sigma),\psi(\hat\sigma))$ be the solution to adjoint state equation \eqref{adjoint} with respect to $\sigma$ and $\hat\sigma$. We denote the sesquilinear form \eqref{se-form} as $\hat{a}(\cdot,\cdot)$ with $\sigma$ replaced by $\hat\sigma$. Then we have
$$
\begin{aligned}
G(\hat\sigma) - G(\sigma) =& \dfrac{1}{2}\mathbf{Re}\{\int_\Gamma(\bE(\hat\sigma) - \bE(\sigma))\cdot n\times\overline{(\bE(\hat\sigma) - \bE^{obs} +\bE(\sigma) - \bE^{obs})}\times n ds\}\\
=-&\dfrac{1}{2}\mathbf{Re}\{\hat{a}(\bF(\hat\sigma),\overline{\bE(\hat\sigma) -\bE(\sigma)}) + b(\nabla\psi(\hat\sigma),\overline{\bE(\hat\sigma)-\bE(\sigma)}) + \\
&\qquad\qquad a(\bF(\sigma),\overline{\bE(\hat\sigma) -\bE(\sigma)})+b(\nabla\psi(\sigma),\overline{\bE(\hat\sigma)-\bE(\sigma)})\}\\
=-&\dfrac{1}{2}\mathbf{Re}\{\hat{a}(\bF(\hat\sigma),\overline{\bE(\hat\sigma) -\bE(\sigma)}) + a(\bF(\sigma),\overline{\bE(\hat\sigma) -\bE(\sigma)})\}\\
=-&\dfrac{1}{2}\mathbf{Re}\{\hat{a}(\bE(\hat\sigma),\overline{\bF(\hat\sigma)})-\hat{a}(\bE(\sigma),\overline{\bF(\hat\sigma)})+a(\bE(\hat\sigma),\overline{\bF(\sigma)})-a(\bE(\sigma),\overline{\bF(\sigma)})\}.
\end{aligned}
$$
By the definition of sesqulinear form (\ref{se-form}),
$$
    \hat{a}(\bE,\bF) =a(\bE,\bF) - \int_{\Omega}i\omega(\hat\sigma - \sigma)\bE\cdot\bF dx
                     = a(\bE,\bF) - it\omega\int_{\Omega_C} \tilde{\sigma}\bE\cdot\bF dx.
$$
Notice $b(\nabla \phi(\hat\sigma) ,\overline{\bF(\hat\sigma)}) =b(\nabla \phi(\sigma) ,\overline{\bF(\hat\sigma)}) =b(\nabla \phi(\hat\sigma) ,\overline{\bF(\sigma)}) =b(\nabla \phi(\sigma) ,\overline{\bF(\sigma)}) = 0$, then by the state equation \eqref{weakeddy}, 
\begin{eqnarray*}
&&\hat{a}(\bE(\hat\sigma),\overline{\bF(\hat\sigma)}) = i\omega\int_{\Omega}\bJ_s \cdot \bF(\hat\sigma)dx.\\
&&\hat{a}(\bE(\sigma),\overline{\bF(\hat\sigma)}) =a(\bE(\sigma),\overline{\bF(\hat\sigma)}) - it\omega\int_{\Omega_C} \tilde\sigma\bE(\sigma)\cdot\bF(\hat\sigma)dx\\
&&\qquad\qquad=i\omega\int_{\Omega}\bJ_s\cdot\bF(\hat\sigma)dx- it\omega\int_{\Omega_C} \tilde\sigma\bE(\sigma)\cdot\bF(\hat\sigma)dx.\\
&&a(\bE(\hat\sigma),\overline{\bF(\sigma)}) = \hat{a}(\bE(\hat\sigma),\overline{\bF(\sigma)})+it\omega\int_{\Omega_C}\tilde\sigma\bE(\hat\sigma)\cdot\bF(\sigma)dx\\
&&\qquad\qquad=i\omega\int_{\Omega}\bJ_s\cdot\bF(\sigma)dx+it\omega\int_{\Omega_C}\tilde\sigma\bE(\hat\sigma)\cdot\bF(\sigma)dx.\\
&&a(\bE(\sigma),\overline{\bF(\sigma)}) = i\omega\int_{\Omega}\bJ_s\cdot \bF(\sigma)dx.
\end{eqnarray*}
Summing up the above equalities, we get
$$G(\hat\sigma)-G(\sigma)=-\frac{1}{2}\mathbf{Re}\{it\omega \int_{\Omega_C}\tilde\sigma \bE(\sigma)\cdot \bF(\hat\sigma)+\tilde\sigma\bE(\hat\sigma)\cdot \bF(\sigma)dx\}.$$
Then
$$
\dfrac{G(\hat\sigma) - G(\sigma)}{t} = \dfrac{1}{2}\omega\mathbf{Im}\{ \int_{\Omega_C}\tilde{\sigma} \bE(\sigma)\cdot\bF(\hat\sigma)dx + \int_{\Omega_C}\tilde{\sigma}\bE(\hat\sigma)\cdot\bF(\sigma)dx\}.
$$
By  $\tilde{\sigma} \in L^{\infty}(\Omega_C)$ and Cauchy inequality,
$$
\begin{aligned}
    |\int_{\Omega}\tilde{\sigma} (\bE(\hat\sigma)-\bE(\sigma))\bF(\sigma)dx |\leq \|\tilde\sigma\|_{L^{\infty}(\Omega_C)}\|\bE(\hat\sigma) - \bE(\sigma)\|_{L^2(\Omega_C)}\|\bF(\sigma)\|_{L^2(\Omega_C)}.
\end{aligned}
$$
We can check that $(\bE(\hat\sigma)-\bE(\sigma),\phi(\hat\sigma)-\phi(\sigma))$ satisfies 
$$
\left\{\begin{aligned}
		a(\bE(\hat\sigma) - \bE(\sigma), \bF) + b(\nabla (\phi(\hat\sigma)-\phi(\sigma)), \bF) dx&=it\omega\int_{\Omega_C}\tilde\sigma \bE(\hat\sigma)\cdot \overline{\bF}dx, \forall \bF \in H_\Gamma(\operatorname{curl};\Omega),\\
		b(\bE(\hat\sigma) - \bE(\sigma) ,\nabla \psi)&=0,  \forall \psi\in H^1_\Gamma(\Omega_0),
\end{aligned}\right.
$$
then by the well-posedness of the above saddle point problem, we have 
$$
\|\bE(\hat\sigma) - \bE(\sigma)\|_{L^2(\Omega)}\leq Ct\omega\|\tilde{\sigma}\bE(\hat\sigma)\|_{L^2(\Omega)}\leq Ct\omega \|\tilde\sigma\|_{L^{\infty}(\Omega_C)}\|\bJ_s\|_{L^2(\Omega)}.
$$
So $\lim\limits_{t\downarrow 0} \|\bE(\hat\sigma) - \bE(\sigma)\|_{L^2(\Omega)} =0 $ and as a conclusion we get
$$
\lim\limits_{t\downarrow 0} \int_{\Omega}\tilde{\sigma} \bE(\hat\sigma)\cdot\bF(\sigma)dx  = \int_{\Omega}\tilde{\sigma} \bE(\sigma)\cdot\bF(\sigma)dx .
$$
Similarly we have 
$$
\lim\limits_{t\downarrow 0} \int_{\Omega}\tilde{\sigma} \bE(\sigma)\cdot\bF(\hat\sigma)dx  = \int_{\Omega}\tilde{\sigma} \bE(\sigma)\cdot\bF(\sigma)dx .
$$
By taking $t\searrow 0$,  we complete the proof with the definition \eqref{gateaux}.
%
\finproof

In this subsection, we will prove the existence and the stability of solution to the optimization problem (\ref{myproblem0}). We also give some important results that will be helpful in convergence analysis in the forthcoming sections. Then we will minimize the objective functional $\bp$ in \eqref{myproblem0} in the convex set $K$. Then we reach the following problem
\begin{equation}\label{myproblem}
	\mathop{\min}_{\sigma \in K}
	\bp(\sigma)= G(\sigma)+\alpha |D\sigma|(\Omega_c) .
\end{equation}

 The following two theorems show the existence and stability of solution to \eqref{myproblem}. The proofs are similar to the corresponding results in \cite{CLZ}, just notice that  $BV(\Omega_C)$ is compactly embedded in $L^1(\Omega_C)$ here while $H^1(\Omega_C)$ is compactly embedded in $L^2(\Omega_C)$ in \cite{CLZ}.  So we omit them here.
\begin{thm}\label{Theorem2.2}
There exists a minimizer $\sigma_\alpha$ to problem (\ref{myproblem}).
\end{thm}

\begin{thm}\label{sta}
Let $\{\bE^k\}$ satisfy $\|\bE^k\times n-\bE^{obs}\times n\|_{L^2(\Gamma)}\to 0$ as $k \to +\infty$. Assume $\sigma_k$ is the minimizer to \eqref{myproblem} with $\bE^{obs}\times n$ replaced by $\bE^k\times n$ in the definition of $G$. Then there is subsequence of $\{\sigma_k\}$ which converges to a minimizer of $\Phi_\alpha(\sigma)$ strongly in $L^1(\Omega_C)$.
\end{thm}
\section{Discrete inverse Problem}\label{sect3}
In this section, we discretize the inverse problem with finite element method and analyze the well-posedness of the discrete problem. 
\subsection{Finite element discretization}
The domain $\Omega$ is partitioned to a family of shape regular, quasi-uniform tetrahedral mesh $\mathcal{M}_h$ and the subscript $h$ denotes the mesh size.  Let $R_h$ be the lowest-order N\'ed\'elec edge element space approximation to $H_\Gamma(\operatorname{curl},\Omega)$. We assume that $\Omega_0$ and $\Omega_C$ are fitted by sub-meshes $\mathcal{M}_0$ and $\mathcal{M}_C$ of $\mathcal{M}_h$, respectively. So we can approximate $H_\Gamma^1(\Omega_0)$ and $H^1_0(\Omega_C)$ by linear Lagrange element on the corresponding sub-meshes and denote the discrete spaces by ${U}_h$ and $V_h$ respectively. 
Let 
$K_h=V_h\cap K$ and 
$D_h$ be the finite element space of piece-wise constant functions on $\Omega_c$. For any $\sigma_h\in K_h$, $\bE_h,\bF_h\in R_h$, we define the discrete sesquilinear form
$$a_h(\bE_h,\bF_h) = \int_{\Omega} \mu^{-1} \nabla \times \bE_h \cdot \nabla \times \overline{\bF_h}-i \omega(\sigma_h+\sigma_0) \bE_h \cdot \overline{\bF_h} d x.$$
Then for given $\sigma_h \in K_h$, the discrete problem to \eqref{weakeddy} is: Find $\bE_h \in R_h, \phi_h \in U_h$, such that
\begin{equation}\label{diseddy}
\left\{\begin{aligned}
		a_h(\bE_h, \bF_h) + b(\nabla\phi_h,\bF_h)&=i \omega \int_{\Omega} \bJ_s \cdot \overline{\bF_h} d x \quad &&\forall \bF_h \in R_h,\\
		b(\bE_h,\nabla\psi_h)&=0 \quad  &&\forall \psi_h\in U_h.
\end{aligned}\right.
\end{equation}
The discrete problem to \eqref{adjoint} is: 
Find $\bF_h \in R_h, \psi_h \in U_h$, such that
\begin{equation}\label{disadjoint}
    \left\{
    \begin{aligned}       a(\bF_h,\tilde{\bE}_h)+b(\nabla\psi_h,\tilde{\bE}_h) =& \int_{\Gamma} \bn\times \overline{(\bE^{obs} - \bE_h)}\times \bn \cdot \overline{\tilde{\bE}_h}ds \qquad &&\forall \tilde{E}_h \in R_h,\\
        b(\bF_h,\nabla\tilde{\phi}_h)=&0 &&\forall \tilde{\phi}_h \in U_h.
    \end{aligned}
    \right.
\end{equation}
Then, for any $\sigma_h\in K_h$, we define the discrete objective functional $G_h:K_h\to R$ as
\begin{equation}
G_h(\sigma_h) = \dfrac{1}{2}\|\bE_h\times n -\bE^{obs}\times n\|_{L^2(\Gamma)}^2, \label{disG}
\end{equation}
where $\bE_h$ is the solution to (\ref{diseddy}) and $\bE^{obs}\times n$ is the data on $\Gamma$. 

 For any element $\mathbf{e}$ in mesh $\mathcal{M}_h$, the number of basis on $\mathbf{e}$ is finite. Hence $\|\cdot\|_{L^p(\mathbf{e})}$ is equivalent to $\|\cdot\|_{L^2(\mathbf{e})}$, i.e. $\|\sigma\|_{L^2(\mathbf{e})}\lesssim\|\sigma\|_{L^p(\mathbf{e})}\lesssim\|\sigma\|_{L^2(\mathbf{e})}$. 
 We list some results from \cite{FEM} here.
 \begin{lem}\label{p2}
 	Assume that the mesh $\mathcal{M}_h$ is shape regular and quasi-uniform. For any $\sigma_h\in V_h$, we have 
  \begin{eqnarray}
\|\sigma_h\|_{L^\infty(\Omega_C)}\leq C_{\infty} h^{-3/2} \|\sigma_h\|_{L^2(\Omega_C)}\\
    \|\sigma_h\|_{L^2(\Omega_C)}\leq C_s h\|\nabla\sigma_h\|_{L^2(\Omega_C)}.\label{scall2}
  \end{eqnarray}
 Where $C_{\infty}$ and $C_s$ are only depended on the shape regularity of mesh $\mathcal{M}_h$ .
 \end{lem}
\debproof
The result can be proved with scaling argument \cite{ciarlet}, we omit the details here. 
\finproof

Similar with the continuous case, we define the G\^ateaux derivative for the discrete objective functional $G_h(\sigma_h)$ as 
\begin{equation}
 (G'_h(\sigma_h),\tilde{\sigma}_h)=\lim_{t\downarrow 0}(G_h(\sigma_h+t\tilde{\sigma}_h)-G_h(\sigma_h))/t.\label{disGh}
\end{equation}
for $\sigma_h,\sigma_h+t\tilde{\sigma}_h\in K_h$. Here and in what follows, $(\cdot, \cdot)$ represents the $L^2$ inner product on $\Omega_C$. Actually, $G'(\sigma_h)$ define the gradient of $G_h$ at $\sigma_h$. 
\begin{lem}\label{disin}
Assume $\sigma_h^1,\sigma_h^2\in K_h$, $\tilde{\sigma}_h \in V_h$ such that $\sigma^1_h+t\tilde{\sigma}_h, \sigma^2_h+t\tilde{\sigma}_h \in K_h$ for some $t>0$, 
 there is a constant $C_{eddy}$ which only depends on the $\bJ_s$, $m,M$ and $\bE^{obs}$,  such that
\begin{eqnarray}
  (G_h'(\sigma_h^1)-G_h'(\sigma_h^2),\tilde{\sigma}_h)\leq C_{eddy}\|\sigma_h^1-\sigma_h^2\|_{L^{\infty}(\Omega_C)}\|\tilde{\sigma}_h\|_{L^\infty(\Omega_C)},\label{disieqK1}\\
G_h(\sigma_h^1)-G_h(\sigma_h^2)-(G_h'(\sigma_h^2),\sigma_h^1-\sigma_h^2)\leq C_{eddy}\|\sigma_h^1-\sigma_h^2\|_{L^\infty(\Omega_C)}^2.\label{disieqK2}
\end{eqnarray}
\end{lem}
 \debproof
 Let $(\bE_h(\sigma_h^i),\phi_h(\sigma_h^i)$  be the solutions to (\ref{diseddy}) with $\sigma_h$ replaced by $\sigma_h^i\in K_h$, $i=1,2$ and  $(\bF_h(\sigma_h^i),\psi_h(\sigma_h^i)$ be the solution to the corresponding adjoint station equation \eqref{disadjoint}.  
Clearly, the estimate
\begin{equation}
\|\bE_h(\sigma^i_h)\|_{H(\operatorname{curl},\Omega)}\leq C\|\bJ_s\|_{L^2(\Omega)}.\label{Ebound}
\end{equation}
is true and the solution to \eqref{disadjoint} has the estimate
\begin{equation*}
\|\bF_h(\sigma^i_h)\|_{H(\operatorname{curl},\Omega)}\leq C\|n\times \overline{(\bE^{obs}-\bE_h)}\times n\|_{L^2(\Gamma)} .
\end{equation*}
Then by trace theorem, 
\begin{eqnarray}
\|\bF_h(\sigma^i_h)\| \leq C(\|\bE_h\|_{H(\operatorname{curl},\Omega)} + \|\bE^{obs}\|_{L^2(\Gamma)})\leq C(\|\bJ_s\|_2+\|\bE^{obs}\|_{L^2(\Gamma)}).\label{Fbound}
\end{eqnarray}
By the same argument in Lemma \ref{derivative}, we have
 \begin{equation}
(G_h'(\sigma_h^i),\tilde{\sigma}_h) = \omega\int_{\Omega_C}\mathbf{Im}(\bE_h(\sigma_h^i)\cdot \bF_h(\sigma_h^i))\tilde{\sigma}_hdx.\label{derivative-Gh}
\end{equation}
Then
\begin{equation}\label{disLip}
	\begin{aligned}
		(G_h'(\sigma_h^1)-G_h'(\sigma_h^2),\tilde{\sigma}_h)
  =&\omega\int_{\Omega_C}\mathbf{Im}(\bE_h(\sigma_h^1)\bF_h(\sigma_h^1)-\bE_h(\sigma_h^2)\bF_h(\sigma_h^2))\tilde{\sigma}_hdx\\
  \leq& \omega |\int_{\Omega_C}(\bE_h(\sigma_h^1)\bF_h(\sigma_h^1)-\bE_h(\sigma_h^2)\bF(\sigma_h^1))\tilde{\sigma}_hdx|+\\
		&\omega |\int_{\Omega_C}(\bE_h(\sigma_h^2)\bF_h(\sigma_h^1)-\bE_h(\sigma_h^2)\bF_h(\sigma_h^2))\tilde{\sigma}_hdx|.
	\end{aligned}    
\end{equation} 
By the H\"older inequality, we have
$$
|\int_{\Omega_C}(\bE_h(\sigma_h^1)\bF_h(\sigma_h^1)-\bE_h(\sigma_h^2)\bF_h(\sigma_h^1))\tilde{\sigma}_hdx|\leq \|\bF_h(\sigma_h^1)\tilde{\sigma}_h\|_{L^2(\Omega_C)}\|\bE_h(\sigma_h^1)-\bE_h(\sigma_h^2)\|_{L^2(\Omega_C)}
$$
Define $a_h^1$ as
$$a_h^1\left(\mathbf{E}_h, \mathbf{F}_h\right)=\int_{\Omega} \mu^{-1} \nabla \times \mathbf{E}_h \cdot \nabla \times \overline{\mathbf{F}_h}-i \omega\left(\sigma_h^1+\sigma_0\right) \mathbf{E}_h \cdot \overline{\mathbf{F}_h} d x \qquad \forall \bE_h,\bF_h\in R_h.$$
Then  $(\bE_h(\sigma_1) - \bE_h(\sigma_2))$ and $(\phi_h(\sigma_1)-\phi_h(\sigma_2))$ satisfy 
$$
\left\{\begin{aligned}
&a_h^1(\bE_h(\sigma_h^1)-\bE_h(\sigma_h^2),\bF_h)+b(\nabla (\phi_h(\sigma_h^1)-\phi_h(\sigma_h^2),\bF_h)\\
&\quad\quad=i\omega\int_{\Omega_C}(\sigma_h^1-\sigma_h^2)\bE_h(\sigma_h^2)\cdot\overline{\bF_h}, \quad \forall \bF_h\in R_h,\\
&b(\bE_h(\sigma_h^1)-\bE_h(\sigma_h^2),\nabla\psi_h)=0, \quad \forall \psi_h\in U_h.
\end{aligned}\right.
$$
Hence by the well-posedness of the above saddle point problem, we can get 
$$
\|\bE_h(\sigma_h^1)-\bE_h(\sigma_h^2)\|_{L^2(\Omega_C)}\leqslant C\|(\sigma_h^1-\sigma_h^2)\bE_h(\sigma_h^2)\|_{L^2(\Omega_C)}.
$$
Then with \eqref{Ebound} and \eqref{Fbound}, we have 
\begin{eqnarray*}
&&	|\int_{\Omega_C}(\bE_h(\sigma_h^1)\bF_h(\sigma_h^1)-\bE_h(\sigma_h^2)\bF_h(\sigma_h^1))\tilde{\sigma}_hdx| \leq \|\bF_h(\sigma_h^1)\tilde{\sigma}_h\|_{L^2(\Omega_C)}\|\bE_h(\sigma_h^1)(\sigma_h^1-\sigma_h^2)\|_{L^2(\Omega_C)}\\		
 &&\qquad \leq \|\bF_h(\sigma_h^1)\|_{L^2(\Omega_C)}\|\tilde{\sigma}_h\|_{L^\infty(\Omega_C)} \|\bE_h(\sigma_h^1)\|_{L^2(\Omega_C)}\|\sigma_h^1-\sigma_h^2\|_{L^\infty(\Omega_C)}\\
&&\qquad	
\leq C\|\mathbf{J}_s\|_{L^2(\Omega_C)}(\|\mathbf{J}_s\|_{L^2(\Omega_C)} +\|\bE^{obs}\|_{L^2(\Gamma)})\|\tilde{\sigma}_h\|_{L^\infty}\|\sigma_h^1-\sigma_h^2\|_{L^\infty(\Omega_C)}.
\end{eqnarray*}
Similarly, for the second term, we can get  
$$
\begin{aligned}
    |\int_{\Omega_C}&(\bE(\sigma_h^2)\bF(\sigma_h^1)-\bE(\sigma_h^2)\bF(\sigma_h^2))\tilde{\sigma}_hdx|\leq\\ &C\|\mathbf{J}_s\|_{L^2(\Omega_C)}(\|\mathbf{J}_s\|_{L^2(\Omega_C)} +\|\bE^{obs}\|_{L^2(\Gamma)})\|\tilde{\sigma}_h\|_{L^\infty}\|\sigma_h^1-\sigma_h^2\|_{L^\infty(\Omega_C)}.
\end{aligned}
$$
This completes the proof of inequality \eqref{disieqK1} by letting $C_{eddy}=C\|\mathbf{J}_s\|_{L^2(\Omega_C)}(\|\mathbf{J}_s\|_{L^2(\Omega_C)} +\|\bE^{obs}\|_{L^2(\Gamma)})$.

Since $\sigma_h^1,\sigma^2_h\in K_h$ and $K_h$ is convex set, we know that  $\sigma^2_h+t(\sigma^1_h-\sigma^2_h)\in K_h$ for $t\in [0,1]$. With the differentiability of $G(\sigma)$, we know that there is some $t\in (0,1)$, such that
$$\begin{aligned}
	G_h(\sigma_h^1)-G_h(\sigma_h^2)-(G_h'(\sigma_h^2),\sigma_h^1-\sigma_h^2)=(G_h'(t\sigma_h^1+(1-t)\sigma_h^2),\sigma_h^1-\sigma_h^2)-(G_h'(\sigma_h^2),\sigma_h^1-\sigma_h^2).
\end{aligned}$$
Then by \eqref{disieqK1}, 
\begin{eqnarray*}
&&G_h(\sigma_h^1)-G_h(\sigma_h^2)-(G_h'(\sigma_h^2),\sigma_h^1-\sigma_h^2)
\leq C_{eddy}\|t\sigma_h^1-t\sigma_h^2\|_{L^\infty(\Omega)}\|\sigma_h^1-\sigma_h^2\|_{L^\infty(\Omega)}\\
&&\qquad \leq C_{eddy}\|\sigma_h^1-\sigma_h^2\|_{L^\infty(\Omega)}^2.
\end{eqnarray*}
This completes the estimate \eqref{disieqK2}.
 \finproof
 
 As a conclusion of Lemma \ref{p2} and Lemma \ref{disin}, we have the following result

\begin{lem}\label{disL2}
Assume $\sigma_h^1,\sigma_h^2\in K_h$, $\tilde{\sigma}_h \in V_h$ such that $\sigma^1_h+t\tilde{\sigma}_h, \sigma^2_h+t\tilde{\sigma}_h \in K_h$ for some $t>0$, 
there is a constant $C_0$ which only depends on the shape regularity of mesh $\mathcal{M}_h$, $\bJ_s$, $K_h$ and $\bE^{obs}$,  such that
\begin{eqnarray}
	(G_h'(\sigma_h^1)-G_h'(\sigma_h^2),\tilde{\sigma}_h)\leq C_0h^{-3}\|\sigma_h^1-\sigma_h^2\|_{L^2(\Omega_C)}\|\tilde{\sigma}_h\|_{L^2(\Omega_C)},\label{newdisieqK1}\\
	G_h(\sigma_h^1)-G_h(\sigma_h^2)-(G_h'(\sigma_h^2),\sigma_h^1-\sigma_h^2)\leq C_0 h^{-3}\|\sigma_h^1-\sigma_h^2\|_{L^2(\Omega_C)}^2.\label{newdisieqK2}
\end{eqnarray}
\end{lem}

\subsection{Well-posedness of discrete inverse problem}
With the finite element discretization, the discrete inverse problem can be recast to the following  optimization problem
\begin{equation}\label{dispro}
	\mathop{\min}_{\sigma_h \in K_h}
	\phi_\alpha(\sigma_h)= G_h(\sigma_h)+\alpha |D\sigma_h|(\Omega_C). 
\end{equation}
The well-posedness of the above problem is given in the following two theorems.

\begin{thm}[existence]\label{disex}
 There exists at least one minimizer to the discrete optimization
problem (\ref{dispro}).
\end{thm}
\debproof
We know $\phi_{\alpha}(\sigma_h)$ is bounded below, porper and lower continuous. Then $\phi_{\alpha}$ has a minimizer  in the closed set $K_h$.
\finproof

\begin{thm}[stability]
Let $\{\bE^k\}$ satisfies $\|\bE^k\times n-\bE^{obs}\times n\|_{L^2(\Gamma)}\to 0$ as $k \to +\infty$. Assume $\sigma_k$ is the minimizer in Theorem \ref{disex}, but with $\bE^{obs}\times n$ replaced by $\bE^k\times n$. Then any limit point of $\{\sigma_n\}$ solves \eqref{dispro}.
\end{thm}
\debproof
 The proof is same as the proof in (\ref{sta}), just notice that $K_h$ is a close set and any limit point of $\{\sigma_h^k\}$ is still in $K_h$.
\finproof

\section{The algorithm and convergence analysis}\label{sect4}
In the section, we will propose a modified ADMM algorithm to solve problem \eqref{dispro} and analyze the convergence. 
In order to deal with the non-smooth and non-convex functional $\phi_\alpha(\sigma_h)$, 
since $|D\sigma|(\Omega_C) = \|\nabla \sigma\|_{L^1(\Omega_C)}$ for $\sigma \in K_h$ since $K_h\subset V_h\subset H^1_0(\Omega_C)\subset W^{1,1}(\Omega_C)$, we recast the problem \eqref{dispro} to the following optimization problem
\begin{equation}\label{dispro1}
\begin{aligned}
		\mathop{\min}_{(\sigma_h,s_h) \in K_h\times K_h}
	\phi_{\alpha,h}(\sigma_h,s_h)=& G_h(\sigma_h)+\alpha \|\nabla s_h\|_{L^1(\Omega_c)},\\
	\text{subject to } s_h =& \sigma_h \text{ in weak sense.}
\end{aligned}
\end{equation}  
To relax the constrain $s_h = \sigma_h$ in \eqref{dispro1}, we introduce the augmented Lagrangian $L_{\alpha,\beta}^h: K_h\times K_h \times(D_h)^3\to \mathbb{R}$ as
\begin{equation}\label{lagran}
\begin{aligned}
	L_{\alpha,\beta}^h(\sigma_h,s_h,y_h)=G_h(\sigma_h)&+\alpha \|\nabla s_h\|_{L^1(\Omega_C)}+(y_h,s_h-\sigma_h)+\frac{\beta}{2} \|\nabla s_h - \nabla\sigma_h\|_{L^2(\Omega_C)}^2.
\end{aligned}
\end{equation} 
Where $\alpha > 0 $ is the regularization parameter and $\beta > 0$ is a penalty parameter.  In the augmented Lagrangian \eqref{lagran},  We  chose $\frac{\beta}{2}\|\nabla s_h - \nabla\sigma_h\|_{L^2(\Omega_C)}^2$ as penalty since this choice can guarantee the existence of minimizers of subproblem \eqref{sigmaiter} in the following algorithm \cite{CLZ}.

Here we remark that if $\sigma$ is approximated with piecewise constant,  $|D\sigma|(\Omega_C) = \|\nabla \sigma\|_{L^1(\Omega_C)}$ is invalid and we can not deal with the total variation regularization term as in \eqref{dispro1}. In that situation, the regularization term should be treated in a different way \cite{chantai}.

Based on the Lagrangian \eqref{lagran}, we introduce a modified ADMM Algorithm \ref{mADMM} to solve problem \eqref{dispro1}.  In each step of the algorithm, we need to solve an optimization sub-problem. 
\begin{algorithm}[htbp]
\caption{Modified ADMM}\label{mADMM}
	\hspace*{0.02in} {\bf Input:} initial values  $\sigma^0,d^0,y^0$\\
	\hspace*{0.02in} {\bf Output:} $\sigma_N$
	\begin{algorithmic}
		\For{$k=0,1,...N-1$} 
		\State 		\begin{equation}\label{sigmaiter}
			\sigma^{k+1}=\mathop{\arg\min}\limits_{\sigma\in K_h} L^h_{\alpha,\beta}(\sigma,s^{k},y^k),
	\end{equation}
 	\begin{equation}\label{s}
			s^{k+1}=\mathop{\arg\min}\limits_{s\in  K_h} L^h_{\alpha,\beta}(\sigma^{k+1},s,y^k),
	\end{equation}
	\begin{equation}\label{y}
		y^{k+1} = G'(\sigma^{k+1})
	\end{equation}
		\EndFor
	\end{algorithmic}
\end{algorithm}


We can solve the first sub-problem by NLCG method \cite{CLZ}.  The optimal condition says
\begin{eqnarray*}
	(\frac{\partial L^h_{\alpha,\beta}}{\partial\sigma}(\sigma^{k+1},s^k,y^k),\tilde\sigma-\sigma^{k+1})\geq 0, \forall \tilde{\sigma}\in K_h.
\end{eqnarray*}
which can be wrote into
\begin{equation}\label{sigmamin}
	(G_h'(\sigma^{k+1}),\widetilde{\sigma}-\sigma^{k+1})\geqslant (y^k, \tilde\sigma - \sigma^{k+1})+(\beta(\nabla s^{k}-\nabla\sigma^{k+1}),\nabla(\widetilde{\sigma}-\sigma^{k+1}))\quad \forall \widetilde{\sigma} \in K_h.
\end{equation}


 
%
Sub-problem \eqref{s} minimizes a strictly convex function on a convex set. There are many algorithms and can be solved rapidly. We use the technique introduced in \cite{ADMM1}, apply an auxiliary variable $d_h \in (D_h)^3$ and change \eqref{s} to equivalent form
\begin{equation}\label{eq2thsub1}
	\mathop{\arg\min}\limits_{d_h \in (D_h)^3, s_h \in K_h,d_h=\nabla s_h} 	\alpha \|d_h\|_{L^1(\Omega_C)} +(y^k,s_h)+ \dfrac{\beta}{2}\|d_h-\nabla\sigma^{k+1} \|_{L^2(\Omega_C)}^2.
\end{equation}
The optimization problem \eqref{eq2thsub1} can be solved by ADMM method, which consists of the iterations
\begin{eqnarray}
	&& d_{t+1} = \mathop{\arg\min}\limits_{d\in (D_h)^3} \alpha \|d\|_{L^1(\Omega_C)} + \dfrac{\beta}{2}\|d - \nabla \sigma^{k+1}\|_{L^2(\Omega_C)}^2 + \dfrac{\rho}{2}\|d - \nabla s_t + u_t\|_{L^2(\Omega_C)}^2,\label{sADMM1}\\ 
	&& s_{t+1} = \mathop{\arg\min}\limits_{s\in K_h} (y^k,s)  + \dfrac{\rho}{2}\|d_{t+1} - \nabla s + u_t\|_{L^2(\Omega_C)}^2,\label{sADMM2}\\
	&& u_{t+1} = u_t + d_{t+1} - \nabla s_{t+1}. \label{sADMM3}
\end{eqnarray}
for given $d_0, s_0, u_0$ . Here $\rho>0$ is a penalty parameter. Convergence of \eqref{sADMM1} - \eqref{sADMM3} is guaranteed by \cite{ADMM1}. The first step \eqref{sADMM1} has a close form solution. The second step \eqref{sADMM2} can be solved by a projection operator. The third step \eqref{sADMM3} can be computed directly. 

The last sub-problem  \eqref{y} in Algorithm \ref{mADMM}  can be computed directly by the adjoint method. It is worth mentioning that the sub-iteration \eqref{y} is different from the traditional ADMM algorithm. Traditionally, the last iteration should be 
\begin{eqnarray}\label{trady}
y^{k+1}=y^k+\beta h,
\end{eqnarray}
where $h\in V_h$ satisfies the following equation
\begin{eqnarray*}
(h,\tilde{\sigma})=(\nabla(s^{k+1}-\sigma^{k+1}),\nabla\tilde{\sigma}) \quad\forall\tilde{\sigma}\in V_h.
\end{eqnarray*}
If the last subproblem is chosen as \eqref{trady}, we fail to prove the convergence of the algorithm. The main difficulty is the inequality optimal condition \eqref{sigmamin} and gradient Lipschitz property of $G'(\sigma)$ can not be used in the convergence proof. So inspired by the convergence analysis we use \eqref{y} as the last sub iteration in Algorithm \ref{mADMM}. We remark that this choice will not bring extra computation cost since solving the first subproblem \eqref{sigmaiter} also need $G'(\sigma^{k+1})$ in the next iteration \cite{CLZ}. 

Before we prove the convergence of Algorithm \ref{mADMM}, we first show some technique inequalities. 
\begin{lem}\label{ieq1}
	$$
	\|\nabla s^k - \nabla \sigma^{k+1}\| _{L^2(\Omega_C)} \leq \beta^{-1} C_0C_sh^{-2}\|\sigma^{k+1} - \sigma^k\|_{L^2(\Omega_C)}
	$$
	$$
	\|s^{k}-\sigma^{k+1}\|_{L^2(\Omega_C)} \leq \beta^{-1}C_0 C_s^2h ^{-1} \|\sigma^{k+1} - \sigma^k\|_{L^2(\Omega_C)}. \label{ineq1}
	$$
\end{lem}
\debproof
Since $s^k$ is in $K_h$. take $\tilde\sigma = s^k$ in \eqref{sigmamin},  use \eqref{newdisieqK1} and \eqref{y},
$$
\begin{aligned}
	\beta\|\nabla s^k-\nabla\sigma^{k+1}\|_{L^2(\Omega_C)}^2\leq & (G'(\sigma^{k+1}) - G'(\sigma^k),s^k-\sigma^{k+1})\\
	\leq& C_0 h ^{-3} \|\sigma^{k+1} - \sigma^k\|_{L^2(\Omega_C)} \| s^k-\sigma^{k+1}\|_{L^2(\Omega_C)}.
\end{aligned}
$$
By Lemma \ref{p2}, 
$$
\|s^k - \sigma^{k+1}\|_{L^2(\Omega_C)}\leq C_s h\|\nabla s^k - \nabla \sigma^{k+1}\|_{L^2(\Omega_C)}
$$
This completes the proof. 

\finproof

Now, we can analyze the decrease of $L^h_{\alpha,\beta}(\sigma^k,d^k,y^k)$. 
\begin{lem}\label{l1}
   There is some positive constant $\beta_0=O(h^{-1})$ such that if  $\beta > \beta_0$, then the series  $\{L_{\alpha,\beta}(\sigma^k,s^k,y^k)\}$ is decreasing. Moreover
\begin{equation}\label{ineq}
	\begin{aligned}
		&L^h_{\alpha,\beta}(\sigma^{k+1},s^{k+1},y^{k+1})-L^h_{\alpha,\beta}(\sigma^k,s^k,y^k)\\
  &\quad\quad \leq -\frac{\beta}{4} C_s^{-2}h^{-2}\|s^{k+1}-s^k\|_{L^2(\Omega_C)}^2
		-(\frac{\beta}{2}C_s^{-2}h^{-2}- C_0 h^{-3} - 2\beta^{-1} C_0^2C_s^2h^{-4})\|\sigma^{k+1}-\sigma^k\|_{L^2(\Omega_C)}^2
	\end{aligned}.
\end{equation}
\end{lem}
\debproof
We have the following decomposition,
\begin{equation}
	\begin{aligned}
		L^h_{\alpha,\beta}(\sigma^{k+1},s^{k+1},y^{k+1})-L^h_{\alpha,\beta}(\sigma^k,s^k,y^k)&=L^h_{\alpha,\beta}(\sigma^{k+1},s^{k+1},y^{k+1})-L^h_{\alpha,\beta}(\sigma^{k+1},s^{k+1},y^k)\\
		&+L^h_{\alpha,\beta}(\sigma^{k+1},s^{k+1},y^{k})-L^h_{\alpha,\beta}(\sigma^{k+1},s^k,y^k)\\
		&+L^h_{\alpha,\beta}(\sigma^{k+1},s^{k},y^{k})-L^h_{\alpha,\beta}(\sigma^k,s^k,y^k)\\
		&=I_1+I_2+I_3.
	\end{aligned}
\end{equation}
Then we estimate $I_i,i=1,2,3$ term by term. For $I_1$, with the help of Lemma \ref{disL2}, Lemma \ref{ineq1} and Young ineqality, for any $\delta >0$, we have
\begin{equation}\label{I1}
	\begin{aligned}
		I_1=&(y^{k+1}-y^k,s^{k+1}-\sigma^{k+1})\\
           =&(G'(\sigma^{k+1} - G'(\sigma^k)) ,s^{k+1} - \sigma^{k+1}) \\
           \leq& C_0h^{-3} \|\sigma^{k+1}-\sigma^{k}\|_{L^2(\Omega_C)} \|s^{k+1}-\sigma^{k+1}\|_{L^2(\Omega_C)} \\
           \leq & C_0h^{-3} \|\sigma^{k+1}-\sigma^{k}\|_{L^2(\Omega_C)}
           (\|s^{k+1}-s^k\|_{L^2(\Omega_C)} + \|s^k-\sigma^{k+1}\|_{L^2(\Omega_C)})\\
           \leq& C_0h^{-3}(\dfrac{\|\sigma^{k+1}-\sigma^k\|_{L^2(\Omega_C)}^2}{4\delta}+ \delta \|s^{k+1} - s^k\|_{L^2(\Omega_C)}^2) + \beta^{-1}C_0^2C_s^2h^{-4} \|\sigma^{k+1} - \sigma^k\|_{L^2(\Omega_C)}^2
	\end{aligned}
\end{equation}
For $I_2$, since $s^{k+1}$ is the solution to \eqref{s}, there is some $\eta^{k+1}$, which is sub-gradient of $g(s)$ at $s^{k+1}$, such that $\forall s\in K_h$,
\begin{equation}\label{smin1}
	\alpha(\eta^{k+1} , s - s^{k+1}) + (y^k, s - s^{k+1}) + \beta (\nabla s^{k+1} - \nabla \sigma^{k+1} , \nabla s - \nabla s^{k+1})\geq 0
\end{equation}
where $g(s) = \|\nabla s\|_{L^1(\Omega_C)}$. Hence we have
$$
\begin{aligned}
	\alpha g(s^k) - \alpha g(s^{k+1}) \geq & \alpha(\eta^{k+1}, s^k - s^{k+1})\\
	\geq &-(y^k, s^k-s^{k+1}) - \beta (\nabla s^{k+1} - \nabla \sigma^{k+1},\nabla s^{k} - \nabla s^{k+1})
\end{aligned}
$$
So we have
\begin{equation}\label{I2}
	\begin{aligned}
		I_2 = &\alpha g(s^{k+1}) - \alpha g(s^k) + (y^k, s^{k+1} - s^k) + \dfrac{\beta}{2} \|\nabla s^{k+1} - \nabla \sigma^{k+1}\|_{L^2(\Omega_C)}^2-\dfrac{\beta}{2}\|\nabla s^k - \nabla \sigma^{k+1}\|_{L^2(\Omega_C)}^2\\
		 \leq & \beta (\nabla s^{k+1} - \nabla \sigma^{k+1},\nabla s^k-\nabla s^{k+1}) +  \dfrac{\beta}{2} \|\nabla s^{k+1} - \nabla \sigma^{k+1}\|_{L^2(\Omega_C)}^2 - \dfrac{\beta}{2}\|\nabla s^k - \nabla \sigma^{k+1}\|_{L^2(\Omega_C)}^2\\
		 =& -\dfrac{\beta}{2} \|\nabla s^{k+1} - \nabla s^k\|_{L^2(\Omega_C)}^2
		 \leq  -\dfrac{\beta}{2} C_s^{-2}h^{-2}\|s^{k+1}-s^k\|_{L^2(\Omega_C)}^2
	\end{aligned}	
\end{equation}
 
Finally, for $I_3$, thanks to Lemma {\ref{disin}}, Lemma \ref{disL2} and \eqref{sigmamin},
\begin{equation}\label{I3}
	\begin{aligned}
		I_3=&G(\sigma^{k+1})-G(\sigma^k)-(y^k,\sigma^{k+1}-\sigma^k)+\dfrac{\beta}{2}\|\nabla s^k-\nabla\sigma^{k+1}\|^2_{L^2(\Omega_C)}-\dfrac{\beta}{2}\|\nabla s^k-\nabla\sigma^k\|_{L^2(\Omega_C)}^2\\
		\leq& (G_h'(\sigma^{k+1}),\sigma^{k+1}-\sigma^k)+C_0h^{-3}\|\sigma^{k+1}-\sigma^k\|_{L^2(\Omega_C)}^2-(y^k,\sigma^{k+1}-\sigma^k)\\
		&-\beta((\nabla s^k-\nabla \sigma ^{k+1}),\nabla\sigma^{k+1}-\nabla\sigma^k)-\frac{\beta}{2}\|\nabla\sigma^{k+1}-\nabla\sigma^k\|_{L^2(\Omega_C)}^2\\
		\leq &-(\frac{\beta}{2}C_s^{-2}h^{-2}-C_0h^{-3})\|\sigma^{k+1}-\sigma^k\|_{L^2(\Omega_C)}^2.\\
	\end{aligned}
\end{equation}
By choosing $\delta = \dfrac{\beta h}{4C_0C_s^2}$ and summing up  \eqref{I1}, \eqref{I2} and \eqref{I3}, we have inequality \eqref{ineq}. The descreasing properties can be confirmed by  choosing $\beta_0$ as 
\begin{equation}
\beta_0=\max\{4C_0C_s^2,2\sqrt{2}C_0C_s^2\}h^{-1}.\label{beta0}
\end{equation}
\finproof

The following lemma tells us that the discrete Lagrangian is bounded below with some conditions.
\begin{lem}\label{l2}
There is some positive constatnt $\beta_1 =O(h^{-1})$, such that $\{L^h_{\alpha,\beta}(\sigma^k,s^k,y^k)\}$ is bounded below for any $\beta >\beta_1$.
\end{lem}
\debproof

By the definition of $L^h_{\alpha,\beta}(\sigma^{k+1},s^{k+1},y^{k+1})$ and \eqref{sigmamin}, then use lemma \ref{p2} and lemma \ref{disL2}, we have
\begin{eqnarray*}
&&		L^h_{\alpha,\beta}(\sigma^{k+1},s^{k+1},y^{k+1})\\
&&\quad=G_h(\sigma^{k+1})+\alpha \|\nabla s^{k+1}\|_{L^1(\Omega_C)}+(y^{k+1}, s^{k+1}- \sigma^{k+1})
		+\frac{\beta}{2}\|\nabla s^{k+1}- \nabla \sigma^{k+1}\|^2_{L^2(\Omega_C)}\\
&&\quad=G_h(\sigma^{k+1})+\alpha \|\nabla s^{k+1}\|_{L^1(\Omega_C)}
+(G_h'(\sigma^{k+1}),s^{k+1}-\sigma^{k+1})
  +\frac{\beta}{2}\|\nabla s^{k+1}-\nabla \sigma^{k+1}\|^2_{L^2(\Omega_C)}\\
&&\quad\geq G_h(\sigma^{k+1})+(G_h'(\sigma^{k+1}), s^{k+1}-\sigma^{k+1})+\frac{\beta C_{s}^{-2}}{2}h^{-2}\| s^{k+1}- \sigma^{k+1}\|^2_{L^2(\Omega_C)}\\
&&\quad\geq  G_h( s^{k+1})-C_0h^{-3}\| s^{k+1} -\sigma^{k+1}\|_{L^2(\Omega_C)}^2
		+\frac{\beta C_{s}^{-2}}{2}h^{-2}\| s^{k+1}- \sigma^{k+1}\|^2_{L^2(\Omega_C)}\\
&&\quad\geq(\dfrac{\beta C_{s}^{-2}}{2}h^{-2}-C_0h^{-3})\|s^{k+1}-\sigma^{k+1}\|_{L^2(\Omega_C)}^2.
\end{eqnarray*}
Then the proof is completed by choosing 
\begin{equation}\label{beta1}
\beta_1 = 2C_0C_{s}^2h^{-1}.
\end{equation}
\finproof

Now we are ready to prove the main result.  

\begin{thm}
If $\beta > \max\{\beta_0,\beta_1\}$, where $\beta_0$ and $\beta_1$  are defined in Lemma \ref{l1} and Lemma \ref{l2} respectively, we have the following convergence results.
\begin{itemize}
	\item  [(a)] $\{L^h_{\alpha,\beta}(\sigma^k,s^k,y^k)\}$ is decreasing, and hence has a limit.
	\item [(b)] We have $\lim\limits_{k\to\infty}\|s^k- \sigma^k\|_{2}\to 0$.
	\item [(c)]If $(\sigma^*,s^*,y^*)$ denotes any limits points of $(\sigma^k,s^k,y^k)$, then
	\begin{equation}\label{convresult}
		\begin{aligned}
			&y^*= G_h'(\sigma^*)\\
			&s^*-\sigma^*=0,\\
			&s^* \in \arg\min\limits_{s\in K_h} \alpha\|\nabla s\|_{L^1(\Omega_C)}+(G_h'(\sigma^*),  s).\\	
		\end{aligned}\\
	\end{equation}
\end{itemize}
\end{thm}
\debproof
	(a)  is obvious from (\ref{ineq}) and Lemma \ref{l2}. \\
	From (\ref{ineq}) and  (\ref{I1}), $\|s^{k+1}-s^k\|_{L^2(\Omega_C)}\to 0,\|\sigma^{k+1}-\sigma^k\|_{L^2(\Omega_C)}\to 0$ as $k\rightarrow\infty$ and hence $\|s^{k}-\sigma^{k+1}\|_{L^2(\Omega_C)} \to 0$, as $k\to \infty$ by Lemma \ref{ineq1}. Then (b) is proved by applying $\|s^{k+1} -\sigma^{k+1}\|_{L^2(\Omega_C)}\leq \|s^k - \sigma^{k+1}\|_{L^2(\Omega_C)} + \|s^{k+1} - s^k\|_{L^2(\Omega_C)}$.\\
	To prove (c), let $\{(\sigma^k,s^k,y^k)\}$ be the sequence generated by Algorithm \ref{mADMM}. Assume its subsequence by $\{(\sigma^{t_k},s^{t_k},y^{t_k})\}$, which converges to $(\sigma^*, s^*, y^*)$ in $L^2(\Omega_C)\times L^2(\Omega_C)\times L^2(\Omega_C)$. Take limit in \eqref{y} we get the first equation of \eqref{convresult}. The second equation of \eqref{convresult} is a direct consequence of (2).		
    As for the third equation, using the optimal condition of \eqref{s}, there is some $\eta^{t_k}$, which is sub-gradient of $g(s)$ at $s^{t_k}$, such that $\forall s\in K_h$,
    \begin{equation}\label{smin2}
    	\alpha(\eta^{t_k} , s - s^{t_k}) + (y^{t_k-1}, s - s^{t_k}) + \beta (\nabla s^{t_k} - \nabla \sigma^{t_k} , \nabla s - \nabla s^{t_k})\geq 0
    \end{equation}
    where $g(s) = \|\nabla s\|_{L^1(\Omega_C)}$. So we have 
    $$
    \begin{aligned}
    	    \alpha g(s) - \alpha g(s^{t_k}) \geq& \alpha(\eta^{t_k} , s - s^{t_k})\\
    	                                    \geq& - (y^{t_k-1}, s- s^{t_k}) - \beta (\nabla s^{t_k} - \nabla \sigma^{t_k} , \nabla s - \nabla s^{t_k})
    \end{aligned}   
    $$
    Take $k\to \infty$, with the help of $\|y^{t_k-1} - y^{t_k}\|_{L^2(\Omega_C)} \lesssim C\|\sigma^{t_k-1} - \sigma^{t_k}\|_{L^2(\Omega_C)}\to 0$ and $\|\nabla s^{t_k-1} - \nabla \sigma^{t_k}\|_{L^2(\Omega)} \lesssim \|\sigma^{t_k-1} - \sigma^{t_k}\|_{L^2(\Omega_C)}\to 0$ in Lemma \ref{ineq1}, we can get
    $$
    \alpha g(s) - \alpha g(s^*) \geq -(y^*, s -s^*)
    $$
 This completes the proof.
\finproof

Here we remark that the penalty parameter $\beta$ can be chosen in $O(h^{-1})$ since $\beta_0,\beta_1$ are both in $O(h^{-1})$.

\section{Numerical experiments}\label{sect5}
In this section, we present some numerical examples to illustrate the efficiency of Algorithm \ref{mADMM}. We set the domain $\Omega=[-2,2]\times [-2,2] \times [-2,0.2]$, with $\Omega_0 = [-2,2]\times [-2,2]\times [0,0.2]$ and $\Omega_c = [-2,2]\times [-2,2]\times [-2,0]$. We implement the algorithms with the parallel hierarchical grid platform (PHG)\cite{PHG}. We carry out all the numerical experiments on a Lenovo notebook with Intel core i5 tenth gen CPU and 8G memory. In our numerical examples, the source is 
$$\mathbf{J}_s=\nabla\times \sum\limits_{i,j=1}^{36} \delta(x-x_{ij}) e_1,$$
where $e_1$ is the unit vector along the $x-axis$ and $x_{ij}= (-2.0 + 0.1\cdot i, -2.0 + 0.1 \cdot j, 0.04)$. Clearly, $\nabla\cdot \mathbf{J}_s = 0$ in weak sense. 
 The algorithms are tested on a mesh $\mathcal{M}_h$ and the data $n\times \bE^{obs}$ is generated by solving problem \eqref{weakeddy} with exact conductivity $\sigma$ on a refined mesh of $\mathcal{M}_h$. We also test our algorithm with noisy data and the noisy data is generated by adding the noise in the following form,
 $$
 n\times \bE^{obs}_{noise} = (1 + \nu \xi) n\times \bE^{obs}.
 $$
 where $\nu$ is noise level and $\xi$ is a uniformly distributed random variable in $[-1,1]$. In the following examples, we choose $\nu = 0.5\%$.
We choose the regularization parameter $\alpha = 10^{-7}$ and the penalty parameter $\beta = 2 \times 10^{-3}$.

In the following examples, the mesh $\mathcal{M}_h$ we used in this example has 745848 edges, which is generated by uniformly refining from an initial mesh with two cubes $\Omega_0$ and $\Omega_C$.  We choose 0 as initial values $\sigma^0=s^0=y^0=0 $ for the algorithm. 
The $\sigma$ sub-problem \eqref{sigmaiter} is solved by NLCG method \cite{CLZ}. To speed up the algorithms, only 3 NLCG iterations is used in Algorithm \ref{mADMM}. For the $s$ sub-problem \eqref{s}, we run the ADMM algorithm \eqref{sADMM1}-\eqref{sADMM3} for 40 iterations, which takes a few computing time since each step can be solved very fast. In the Algorithm \ref{mADMM}, we choose the truncation parameter $m= 0.05\times (k-1)$ and $M=15$ in the $k$-th iteration.

\subsection{Example 1}
In this example, the background conductivity $\sigma_0 = 1$ in $\Omega_C$ and the object $\Omega_2= [-0.3,0.3]\times [-0.3,0.3]\times[-1.0,-0.4]$ . The exact abnormal conductivity is given by $\sigma =5$ in $\Omega_2$ and $\sigma=0$ in $\Omega_C\backslash\bar\Omega_2$. So the exact conductivity is 6 in $\Omega_2$ and 1 in $\Omega_C\backslash\Omega_2$.

\begin{figure}[htbp]
	\centering
	\subfigure
	{
		\begin{minipage}{0.45\linewidth}
			\centering
			\includegraphics[scale=0.2]{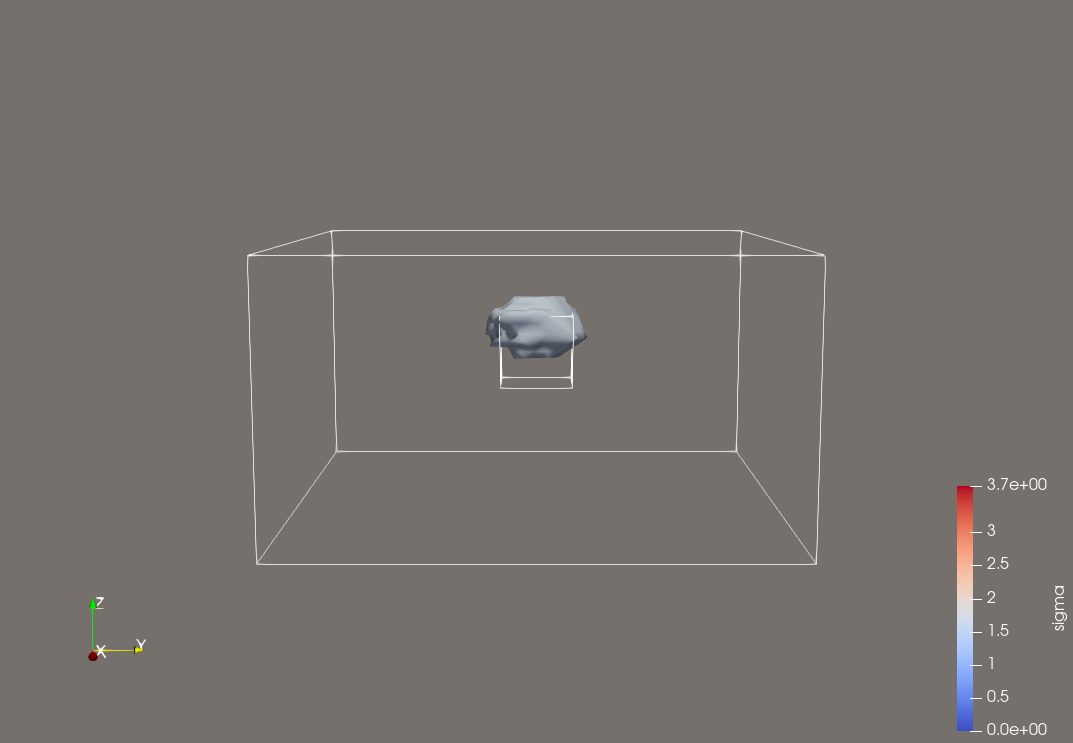}
		\end{minipage}
	}	
	\subfigure
	{
		\begin{minipage}{0.45\linewidth}
			\centering
			\includegraphics[scale=0.2]{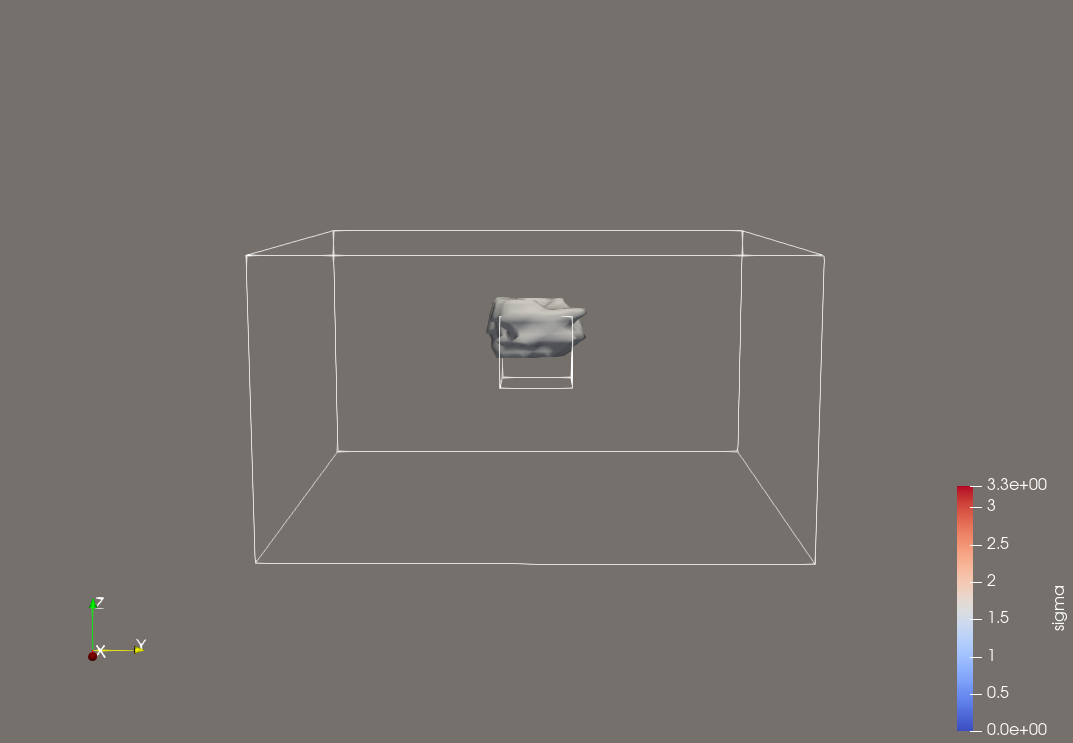}
		\end{minipage}
	}
	\caption{The results for example 1 without noise (left) and with noise (right). }\label{E1}
\end{figure}

\begin{figure}[htbp]
	\centering
     \includegraphics[scale=0.6]{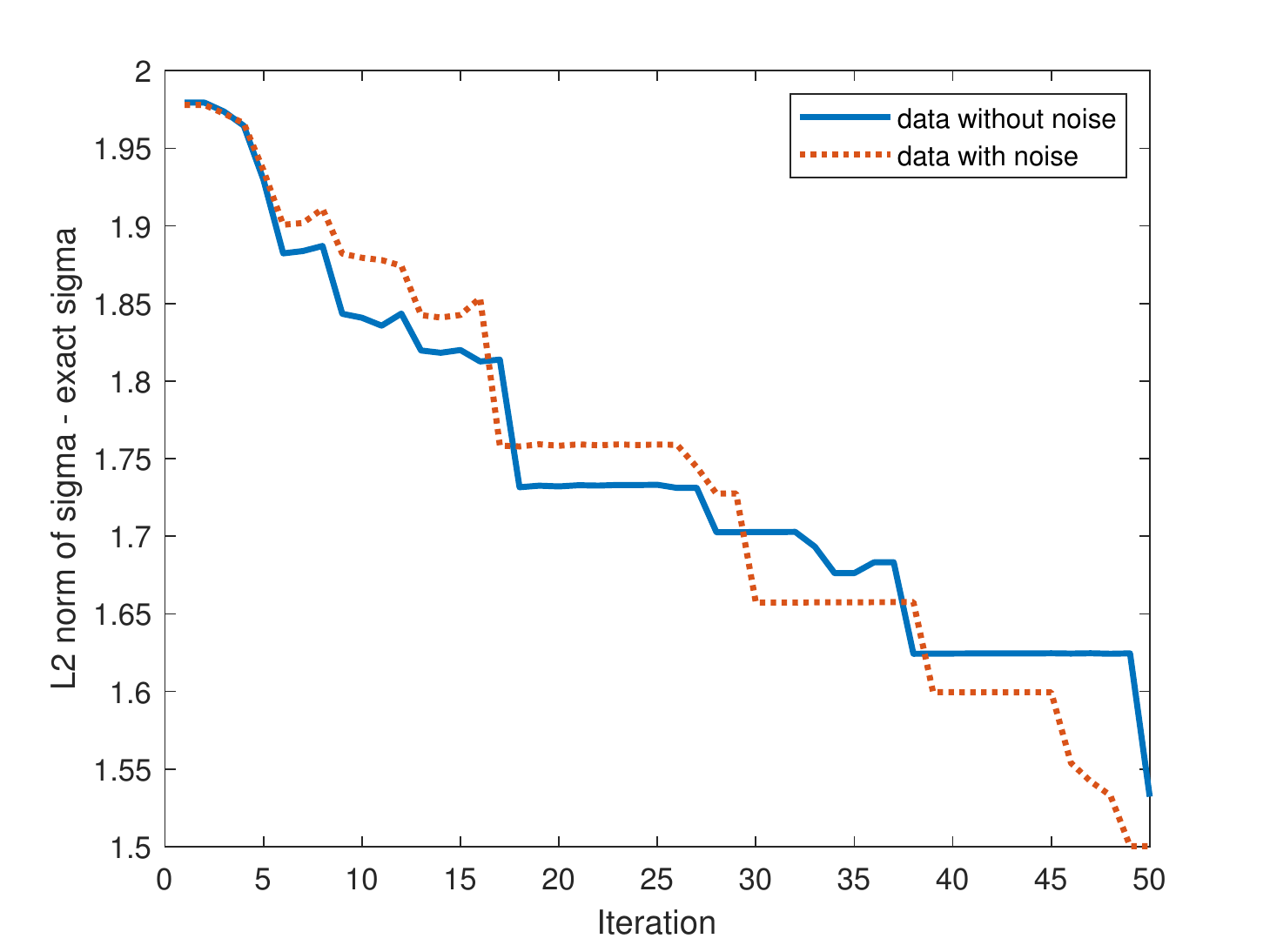}
    \caption{The decreasing of $L_2$ errors for the example 1.}\label{E1error}
\end{figure}

Figure \ref{E1} shows iso-surface of the recovered $\sigma$ with iso-value 1.6 in 50 iterations by Algorithm \ref{mADMM}.  The $L^2$ errors of $\sigma - \sigma^{real}$ with respect to the iterations are shown in Figure \ref{E1error}.  The results show that the algorithm can recover the domain and value simultaneously then performs well  with and without noise.

\subsection{Example 2}
\begin{figure}[htbp]
	\centering
	\subfigure
	{
		\begin{minipage}{0.45\linewidth}
			\centering
			\includegraphics[scale=0.2]{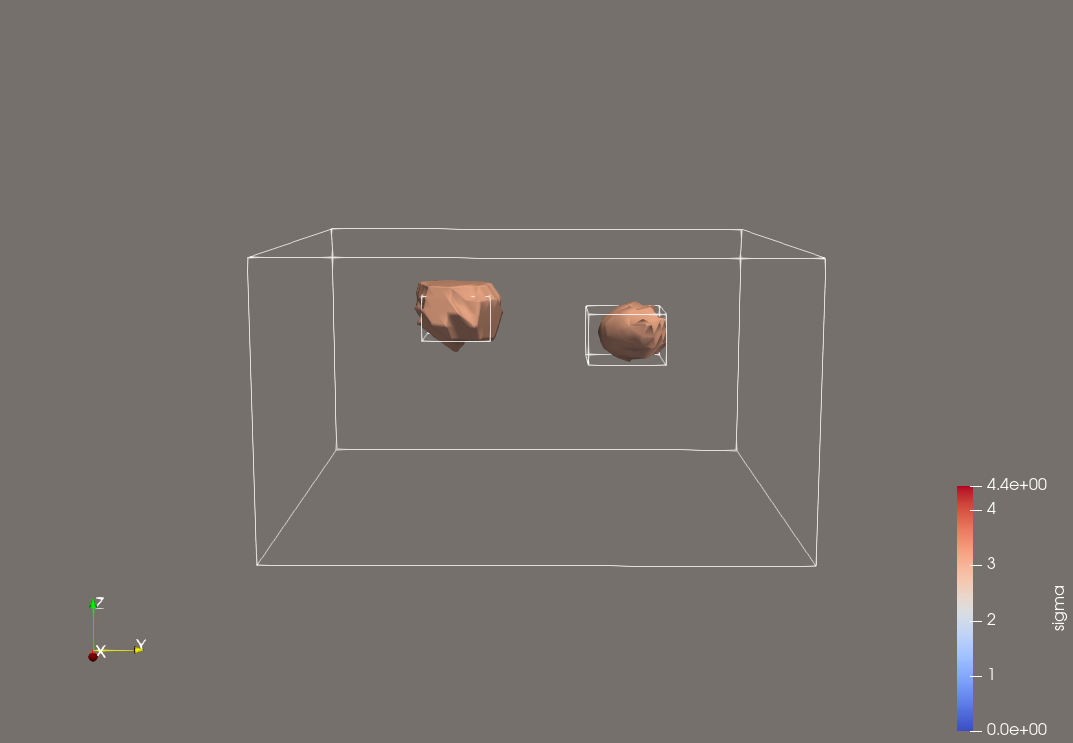}
		\end{minipage}
	}	
	\subfigure
	{
		\begin{minipage}{0.45\linewidth}
			\centering
			\includegraphics[scale=0.2]{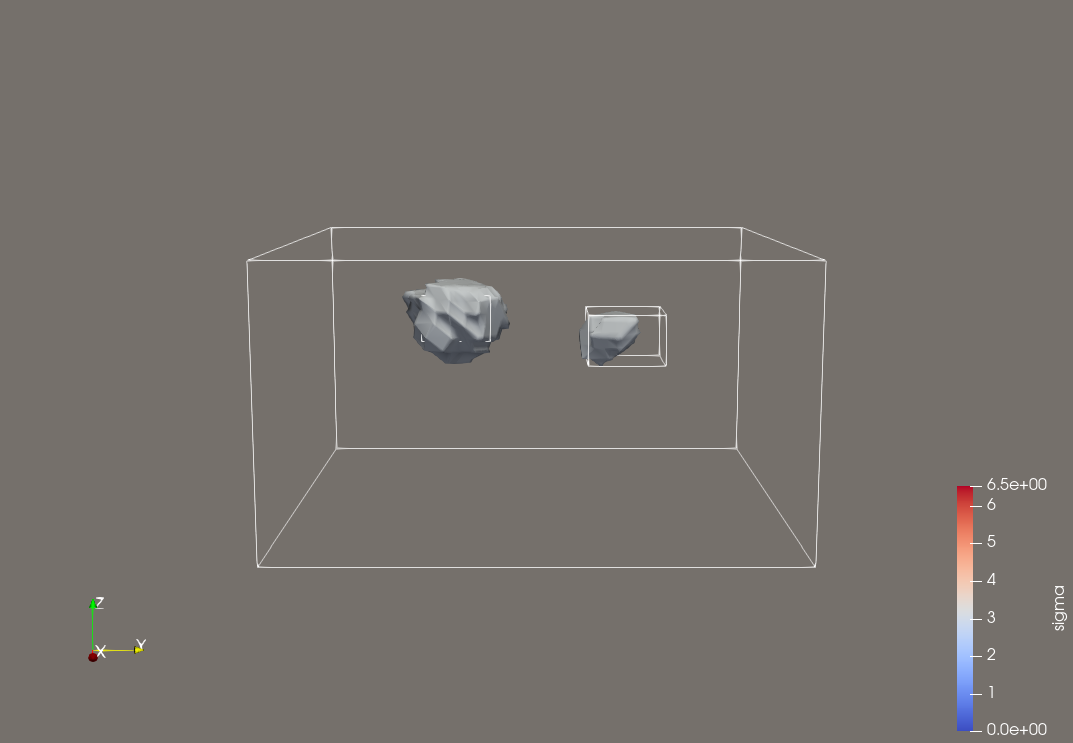}
		\end{minipage}
	}
	\caption{The results for example 2 without noise (left) and with noise (right). }\label{E2}
\end{figure}

\begin{figure}[htbp]
	\centering
	\includegraphics[scale=0.6]{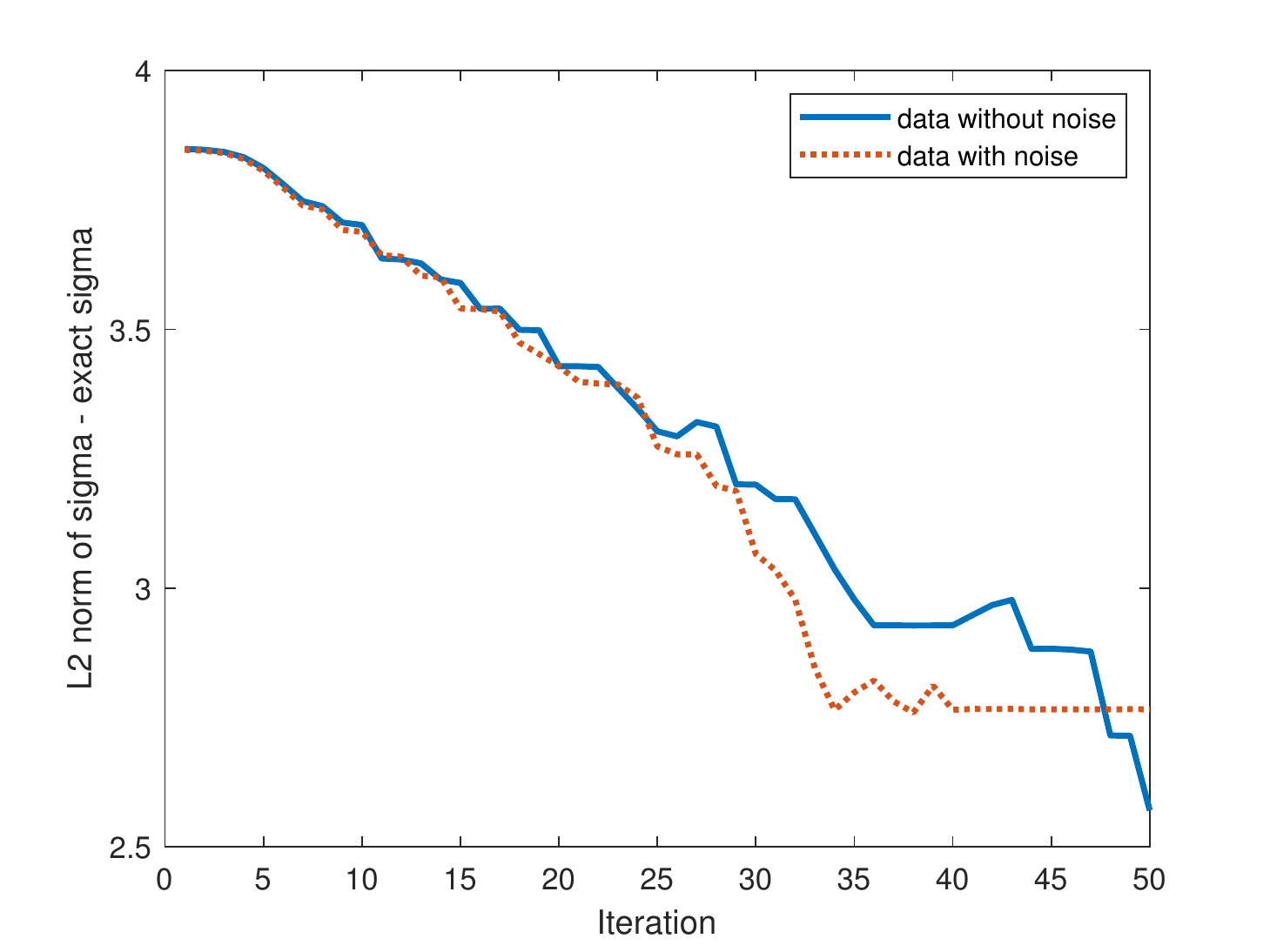}
	\caption{The decreasing of $L_2$ errors for the example 2.}\label{E2error}
\end{figure}

In this example, the background conductivity $\sigma_0 = 1$ in $\Omega_C$ and the abnormal domain $\Omega_2$ consists of two subdomains $D_1= [-0.1,-0.4]\times [-1.0, -0.4]\times[-0.7,-0.3] $ and $D_2=[0.4,1.0] \times [0.4,1.0] \times [-0.7, -0.3]$. 
In this example, we set the abnormal conductivity as follows.
$$
\sigma = \left\{
\begin{aligned}
	10& \qquad \text { in } D_1\\
	6& \qquad \text{ in } D_2 \\
	0& \qquad \text { in } \Omega_C\backslash\Omega_2
\end{aligned}\right.
$$
Figure \ref{E2} shows iso-surface of the recovered $\sigma$ with iso-value 3.0 in 50 iterations by Algorithm \ref{mADMM}. The left picture is the recovery result without noise and the right is the result with noise. In each picture, we can find that the recovered domain at $D_1$ is bigger than that at $D_2$, Considering that we show the recoveries with same iso-value, this means that the recovered value is bigger in $D_1$ and smaller in $D_2$.  Then Figure \ref{E2}  tells us that the algorithm can recover the different conductivity in separative domains. The $L^2$ errors of $\sigma - \sigma^{real}$ with clear data and noisy data are showed in Figure \ref{E2error}. Just as in Example 1, our algorithm exhibits good convergence performance in this example.

\subsection{Example 3}

In this example, we test the situation that the abnormal domain is more complicated. We choose $\Omega_2$ to be a L-shape domain. The background conductivity $\sigma_0 = 1$ in $\Omega_C$. We let $\Omega_2 = [-1.5,-1.0]\times[-1.5,1.5]\times[-1.0,0.4]\cup[-1.0,1.5]\times[-1.5,-1.0]\times[-1.0,-0.4]$. The exact abnormal conductivity is given by $\sigma =10$ in $\Omega_2$ and $\sigma=0$ in $\Omega_C\backslash\bar\Omega_2$.

\begin{figure}[htbp]
	\centering
	\subfigure
	{
		\begin{minipage}{0.45\linewidth}
			\centering
			\includegraphics[scale=0.2]{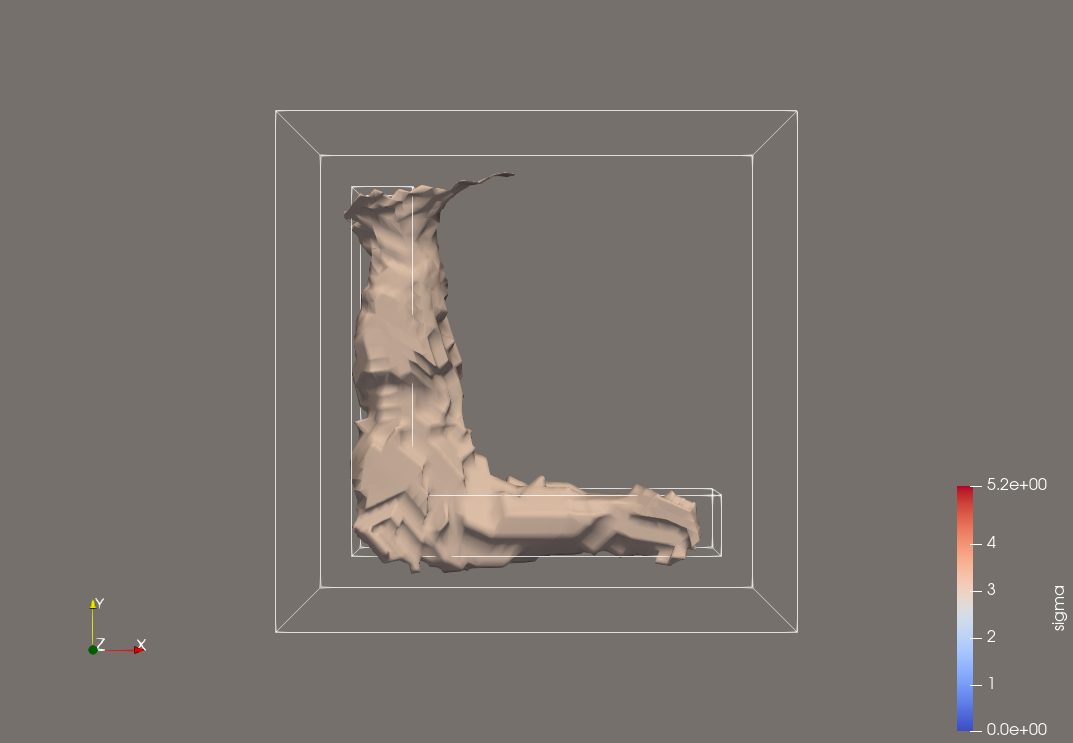}
		\end{minipage}
	}	
	\subfigure
	{
		\begin{minipage}{0.45\linewidth}
			\centering
			\includegraphics[scale=0.2]{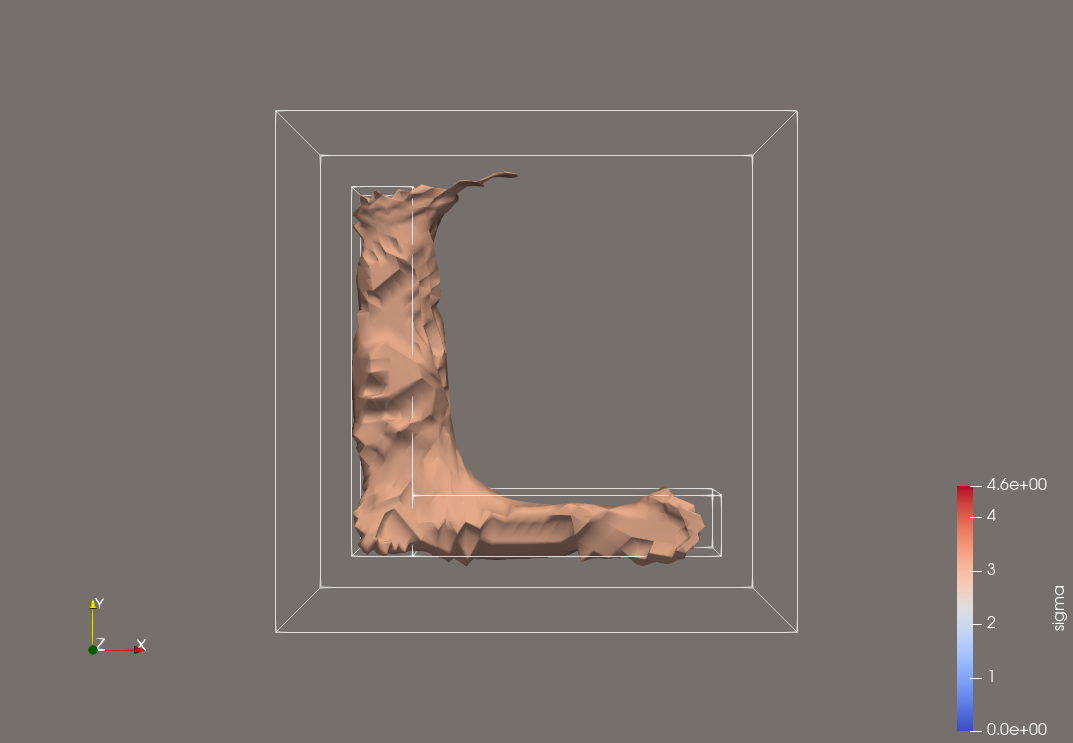}
		\end{minipage}
	}
	\caption{The results for example 3 without noise (left) and with noise (right). }\label{E3}
\end{figure}

\begin{figure}[htbp]
	\centering
	\includegraphics[scale=0.6]{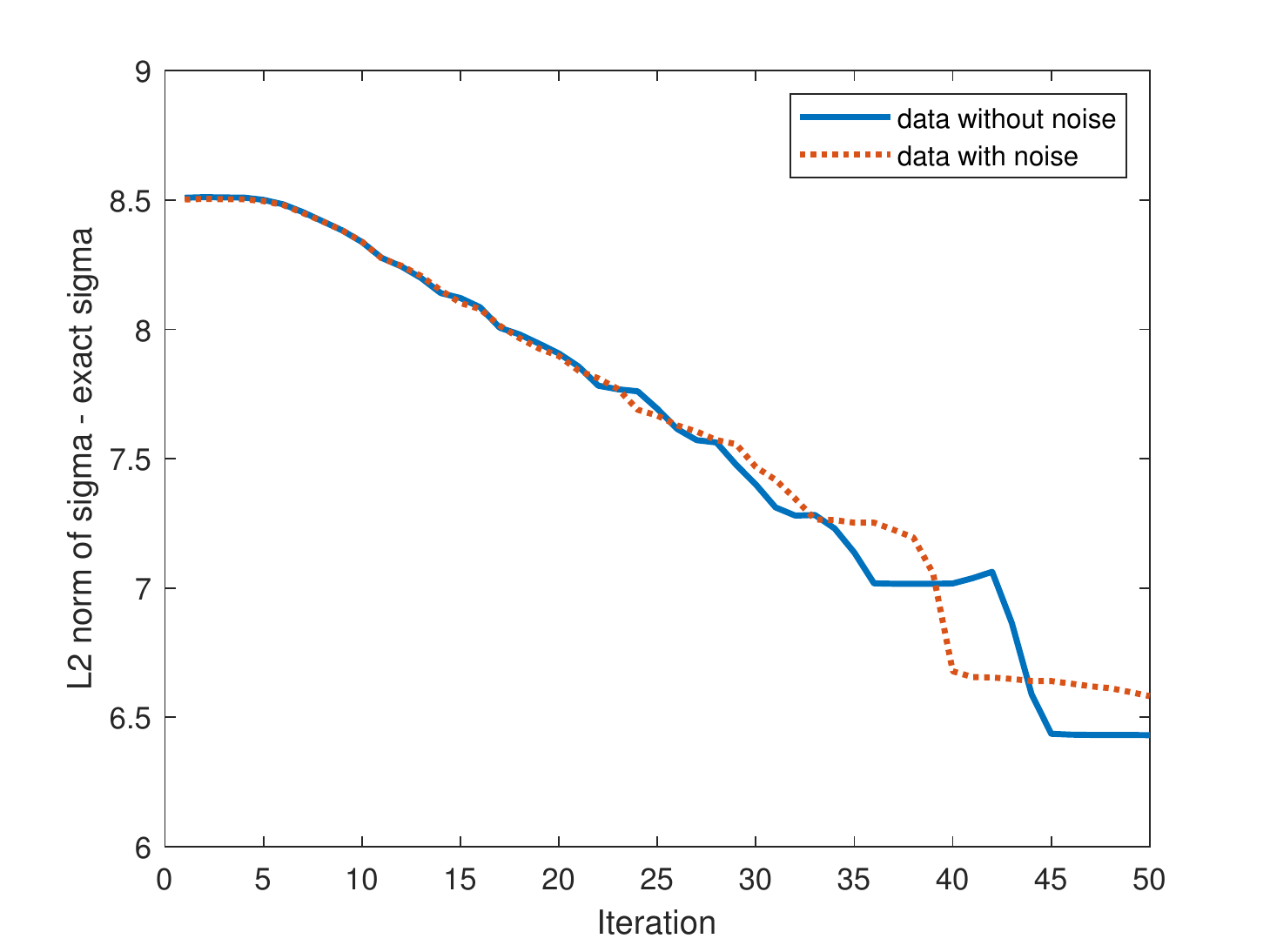}
	\caption{The decreasing of $L_2$ errors for the example 3.}\label{E3error}
\end{figure}

 Figure \ref{E3} shows iso-surface of the recovered $\sigma$ with iso-value 3.0 in 50 iterations by Algorithm \ref{mADMM}. The $L^2$ errors of $\sigma - \sigma^{real}$ with respect to iteration for data with and without noise are showed in Figure \ref{E3error}.  The results show that our algorithm can recover the L-shaped conductor and converge very well, too.

\section{Concluding remarks}\label{sect6}
We study the inverse eddy current problem by assuming that the conductivity lies in bound variation space. Firstly, we analyze the well-posedness of the inverse problem with total variation regularization. Then  we discrete the inverse problem with finite element method and give the well-posedness of the discrete inverse problem. In the framework of ADMM, we propose Algorithm \ref{mADMM} to solve the corresponding inverse problem.  We emphasise that our algorithm is different with the traditional ADMM algorithm since the last sub-problem \eqref{y} is quite different.  Then we can prove the convergence by virtue of the  gradient Lipschitz property of the discrete objective function.  Finally, we give some examples to illustrate the efficiency of the proposed algorithm. 

 \end{document}